\documentclass[fleqn]{article}
\usepackage{amssymb,latexsym,theorem}
\newtheorem{theorem}{Theorem}[section]
\newtheorem{proposition}[theorem]{Proposition}
\newtheorem{lemma}[theorem]{Lemma}
\newtheorem{corollary}[theorem]{Corollary}
\newtheorem{example}[theorem]{Example}

\title{Embeddings of Line-grassmannians of Polar Spaces in Grassmann Varieties}
\author{A. Pasini and I. Cardinali}
\date{}

\begin{document}
\maketitle
\begin{abstract}
An embedding of a point-line geometry $\Gamma$ is usually defined as an injective mapping $\varepsilon$ from the
point-set of $\Gamma$ to the set of points of a projective space such that $\varepsilon(l)$ is a projective line
for every line $l$ of $\Gamma$, but different situations have lately been considered in the literature, where
$\varepsilon(l)$ is allowed to be a subline of a projective line or a curve. In this paper we propose a more
general definition of embedding which includes all the above situations and we focus on a class of embeddings,
which we call Grassmman embeddings, where the points of $\Gamma$ are firstly associated to lines of a projective
geometry $\mathrm{PG}(V)$, next they are mapped onto points of $\mathrm{PG}(V\wedge V)$ via the usual projective
embedding of the line-grassmannian of $\mathrm{PG}(V)$ in $\mathrm{PG}(V\wedge V)$. In the central part of our
paper we study sets of points of $\mathrm{PG}(V\wedge V)$ corresponding to lines of $\mathrm{PG}(V)$ totally
singular for a given pesudoquadratic form of $V$. Finally, we apply the results obtained in that part to the
investigation of Grassmann embeddings of several generalized quadrangles.
\end{abstract}

\section{Introduction}

According to a well established definition, an embedding of a point-line geometry $\Gamma$ is an injective
mapping $\varepsilon$ from the point-set $\cal P$ of $\Gamma$ to the set of points of the projective geometry
$\mathrm{PG}(V)$ of subspaces of $V$ such that $\varepsilon({\cal P})$ spans $\mathrm{PG}(V)$ and
$\varepsilon(l)$ is a projective line, for every line $l$ of $\Gamma$. However, different interesting situations
have lately been considered, where $\varepsilon(l)$ is allowed to be a subline of a projective line or a curve
(a conic or a unital, for instance) or even a surface. In Section \ref{Intro:Emb} we propose a more general
definition of embedding which includes all the above situations and we sketch a bit of theory for it. More
general definitions can be considered, where the embedding is hosted by a group instead of a projective space
(see \cite{PasEE}, for instance), but in this paper we are not willing to go so far. At the end of Section
\ref{Intro:Emb} we consider a class of embeddings, which we call Grassmman embeddings, where the points of
$\Gamma$ are firstly associated to lines of $\mathrm{PG}(V)$ for a vector space $V$ of finite dimension $n$ over
a commutative division ring $\mathbb{F}$, next they are mapped onto points of $\mathrm{PG}(V\wedge V)$ via the
usual projective embedding of the line-grassmannian of $\mathrm{PG}(V)$ in $\mathrm{PG}(V\wedge V)$.

In Section \ref{Section 3} we study sets of points of $\mathrm{PG}(V\wedge V)$ corresponding to lines of
$\mathrm{PG}(V)$ totally isotropic or totally singular for a given sesquilinear or quadratic form. In
particular, we describe those sets as solutions of suitable sets of equations. We do that with the help of the
following trivial observation: $V\wedge V$ is isomorphic to the vector space of anti-symmetric matrices of order
$n$ with entries in $\mathbb{F}$. In spite of its triviality, this remark is extremely useful to write down
equations for the sets we are interested in. Indeed it allows to pack a number of scalar equations in one single
matrix equation by which equations for tangent spaces can be deduced quite easily. We do not know if this trick
has ever been used before. Perhaps it has, but we have not found any track of it in the literature we are aware
of.

In the last section of this paper we apply the results of Section \ref{Section 3} to Grassmann embeddings of a
number of generalized quadrangles. We also compare those embeddings with other embeddings, which we call
veronesean embeddings, obtained by composing a projective embedding $\varepsilon$ with the usual quadric
veronesean mapping of the projective space hosting $\varepsilon$.

To finish this introduction we fix some conventions to be used throughout this paper.

\bigskip

\noindent {\bf Notation.} We will often go back and forth from a vector space $V$ to its projective geometry
$\mathrm{PG}(V)$. In order to avoid any confusion, we keep the usual symbol $\langle . \rangle$ for spans in $V$
while we denote spans in $\mathrm{PG}(V)$ by the symbol $\langle . \rangle^{\mathrm{pr}}$. Given a non-zero
vector $v\in V$ we denote by $[v]$ the point of $\mathrm{PG}(V)$ represented by $v$. Given a set $X$ of vectors
of $V$ we put $[X] = \{[v]\}_{v\in X, v\neq 0}$. In particular, if $X$ is a subspace of $V$ then $[X]$ is a
subspace of $\mathrm{PG}(V)$ (and $\mathrm{dim}(X) = \mathrm{dim}([X])+1$). On the other hand, let $X$ be a set
of points of $\mathrm{PG}(V)$. Then $X$ is a set of $1$-dimensional subspaces of $V$. By a harmless abuse we can
switch from that set of subspaces to its union $\cup_{p\in X}p$, thus regarding $X$ as a subset of $V$.
Accordingly, we write $\langle X\rangle$ for $\langle \cup_{p\in X}p\rangle$.

\section{A generalized notion of embedding}\label{Intro:Emb}

Throughout this section $\Gamma = ({\cal P},{\cal L})$ is a point-line geometry, $\cal P$ is its set of points
and $\cal L$ its set of lines. The members of $\cal L$ are regarded as subsets of $\cal P$, as usual. We assume
that $\Gamma$ is connected and no two distinct lines of $\Gamma$ meet in more than one point.

\subsection{Definitions}\label{Basics-proj}

Given a point-line geometry $\Gamma = ({\cal P},{\cal L})$, a vector space $V$ and a positive integer $d$, a
{\em locally $d$-dimensional embedding} ($d$-{\em embedding} for short) of $\Gamma$ in $\mathrm{PG}(V)$ is an
injective mapping $\varepsilon$ from $\cal P$ to the set of points of $\mathrm{PG}(V)$ satisfying the following:

\begin{itemize}
\item[$(E1)$] for every line $l\in{\cal L}$ the image $\varepsilon(l) := \{\varepsilon(p)\}_{p\in l}$ of $l$
spans a $d$-dimensional subspace $\langle \varepsilon(l)\rangle^{\mathrm{pr}}$ of $\mathrm{PG}(V)$.
\item[$(E2)$] $\langle \varepsilon(l)\rangle^{\mathrm{pr}}\cap \varepsilon({\cal P}) = \varepsilon(l)$ for every
line $l\in{\cal L}$; \item[$(E3)$] $\langle \varepsilon({\cal P})\rangle^{\mathrm{pr}} = \mathrm{PG}(V)$.
\end{itemize}

Let $l$ and $m$ be distinct lines of $\Gamma$. The injectivity of $\varepsilon$ implies that
$|\varepsilon(l)\cap\varepsilon(m)| \leq 1$ while property $(E2)$ implies that $\langle
\varepsilon(l)\rangle^{\mathrm{pr}} \neq \langle \varepsilon(m)\rangle^{\mathrm{pr}}$.

Note also that, by $(E1)$, if $\Gamma$ admits a $d$-embedding then every line of $\Gamma$ has at least $d+1$ points.

We write $\varepsilon:\Gamma\stackrel{d}{\rightarrow}\mathrm{PG}(V)$ to mean that $\varepsilon$ is a
$d$-embedding of $\Gamma$ in $\mathrm{PG}(V)$. We call $d$ the {\em local dimension} of $\varepsilon$. The {\em
global dimension} $\mathrm{dim}(\varepsilon)$ of $\varepsilon$ (also {\em dimension} of $\varepsilon$ for short)
is the dimension of $\mathrm{PG}(V)$, but we warn that many authors, including ourselves in earlier papers, when
dealing with particular classes of $d$-embeddings as projective or quadratic embeddings (to be defined below),
take $\mathrm{dim}(V)$ instead of $\mathrm{dim}(\mathrm{PG}(V))$ as the dimension of $\varepsilon$. If
$\mathbb{F}$ is the underlying division ring of $V$ we say that $\varepsilon$ is {\em defined over}
$\mathbb{F}$, also that $\varepsilon$ is an $\mathbb{F}$-{\em embedding of local dimension} $d$, or an
$\mathbb{F}$-$d$-{\em embedding}, for short.

Let $\varepsilon$ be a $1$-embedding. Then $\varepsilon$ maps the lines of $\Gamma$ into lines of
$\mathrm{PG}(V)$. We say that $\varepsilon$ is {\em full} if $\varepsilon(l) = \langle
\varepsilon(l)\rangle^{\mathrm{pr}}$ for every line $l\in{\cal L}$. If $\varepsilon$ is not full then we say
that it is {\em lax}.

Full $1$-embeddings are often called {\em projective embeddings} in the literature. We shall follow this habit
in this paper, but we warn the reader that many authors use a different terminology, calling $1$-embeddings `lax
projective embeddings' (also just lax embeddings, as in Van Maldeghem \cite{HVM}) and projective embeddings in
our sense `full projective embeddings'.

We also slightly change our notation as follows: if $\varepsilon$ is a projective embedding then we simply write
$\varepsilon:\Gamma\rightarrow\mathrm{PG}(V)$, without keeping a record of the local dimension $d = 1$ in our
notation.

A {\em quadratic embedding} is a $2$-embedding $\varepsilon$ defined over a field (namely a commutative division
ring) and such that $\varepsilon(l)$ is a non-singular conic of the projective plane $\langle
\varepsilon(l)\rangle^{\mathrm{pr}}$, for every line $l$ of $\Gamma$. Quadratic embeddings are called veronesean
embeddings in \cite{CP1} and \cite{CP2}.

Not all $2$-embeddings are quadratic. For instance, a projective space admits many non-quadratic $2$-embeddings
(see \cite{TVM} and \cite{TVMmore}, where these embeddings are called generalized veronesean and lax generalized
veronesean embeddings). The $2$-embeddings of $\mathrm{PG}(n,\mathbb{F})$ in $\mathrm{PG}(m,\mathbb{K})$ for a
division ring $\mathbb{K}$ possibly different from $\mathbb{F}$, with $|\mathbb{F}| > 2$ and $m$ at least
${{n+2}\choose 2}-1$ (which is the dimension of the usual quadric veronesean embedding of
$\mathrm{PG}(n,\mathbb{F})$) are classified in \cite{TVM} and \cite{TVMmore}. In particular, it is proved that
$\mathbb{F}\subseteq \mathbb{K}$ and $m = {{n+2}\choose 2}-1$ in any case.

\subsubsection{Morphisms and quotients}\label{morphisms}

Morphisms, quotients and hulls can be defined for arbitrary $d$-embeddings just in the same way as for
projective embeddings. We consider morphisms and quotients in this subsection. In the next one we shall turn to
hulls.

Given two $\mathbb{F}$-embeddings $\varepsilon_1:\Gamma\stackrel{d_1}{\rightarrow}\mathrm{PG}(V_1)$ and
$\varepsilon_2:\Gamma\stackrel{d_2}{\rightarrow}\mathrm{PG}(V_2)$, a {\em morphism} $f:\varepsilon_1\rightarrow
\varepsilon_2$ from $\varepsilon_1$ to $\varepsilon_2$ is a semi-linear mapping $f:V_1\rightarrow V_2$ such that
$\varepsilon_2 = f\cdot\varepsilon_1$. To be precise, in this definition we should refer to the mapping
$\mathrm{PG}(f)$ from the attenuated space $\mathrm{PG}(V_1)\setminus[\mathrm{Ker}(f)]$ to $\mathrm{PG}(V_2)$
induced by $f$ rather than to $f$ itself (compare Kasikova and Shult \cite{KS}), but we prefer to take $f$ in
place of $\mathrm{PG}(f)$. This is an abuse, but it is harmless. It makes a few statements slightly clumsy (see
for instance the description of $\mathrm{Aut}(\varepsilon)$, a few lines below), but it saves us the trouble of
introducing attenuated spaces.

Note that, since $\langle \varepsilon_2({\cal P})\rangle^{\mathrm{pr}} = \mathrm{PG}(V_2)$, the equality
$\varepsilon_2 = f\varepsilon_1$ forces the mapping $f:V_1\rightarrow V_2$ to be surjective. If $f$ is bijective
then we say that $f$ is an {\em isomorphism} from $\varepsilon_1$ to $\varepsilon_2$. If $\varepsilon_1$ and
$\varepsilon_2$ are isomorphic in this sense, we write $\varepsilon_1\cong \varepsilon_2$. If a morphism exists
from $\varepsilon_1$ to $\varepsilon_2$ then we write $\varepsilon_1\geq \varepsilon_2$ and we say that
$\varepsilon_2$ is a {\em morphic image} of $\varepsilon_1$. If moreover $\varepsilon_1\not\cong\varepsilon_2$
then we say that $\varepsilon_2$ is {\em proper morphic image} of $\varepsilon_1$ and we write $\varepsilon_1 >
\varepsilon_2$.

The automorphisms of a $d$-embedding $\varepsilon:\Gamma\stackrel{d}{\rightarrow}\mathrm{PG}(V)$ form a group
$\mathrm{Aut}(\varepsilon)$, equal to the point-wise stabilizer of $\varepsilon({\cal P})$ in the group
$\Gamma\mathrm{L}(V)$ of all invertible semilinear transformations of $V$. When $\varepsilon$ is projective or
quadratic then $\mathrm{Aut}(\varepsilon)$ is equal to the center $Z(\mathrm{GL}(V))$ of $\mathrm{GL}(V)$ (see
\cite[Proposition 9]{PasVM} for a proof of this claim when $\varepsilon$ is projective), but in general
$\mathrm{Aut}(\varepsilon)\neq Z(\mathrm{GL}(V))$. Given two $d$-embeddings
$\varepsilon_1:\Gamma\stackrel{d}{\rightarrow}\mathrm{PG}(V_1)$ and
$\varepsilon_2:\Gamma\stackrel{d}{\rightarrow}\mathrm{PG}(V_2)$ such that $e_1\geq e_2$, let
$f:\varepsilon_1\rightarrow\varepsilon_2$ be a morphism. Then $f \cdot \mathrm{Aut}(\varepsilon_1)$ is the
family of all morphisms from $\varepsilon_1$ to $\varepsilon_2$. (Note that $\mathrm{Aut}(\varepsilon_2)\cdot f
= f\cdot \mathrm{Aut}(\varepsilon_1)_{\mathrm{Ker}(f)}$ where $\mathrm{Aut}(\varepsilon_1)_{\mathrm{Ker}(f)}$ is
the setwise stabilizer of $\mathrm{Ker}(f)$ in $\mathrm{Aut}(\varepsilon_1)$, whence
$\mathrm{Aut}(\varepsilon_2)\cdot f\cdot \mathrm{Aut}(\varepsilon_1) = f\cdot\mathrm{Aut}(\varepsilon_1)$.) In
particular, when $\varepsilon_1$ and $\varepsilon_2$ are projective or quadratic and $\varepsilon_1\geq
\varepsilon_2$, the morphism $f:\varepsilon_1\rightarrow\varepsilon_2$ is unique modulo scalars.

Given a $d$-embedding $\varepsilon:\Gamma\stackrel{d}{\rightarrow}\mathrm{PG}(V)$, let $K$ be a
subspace of $V$ satisfying the following:

\begin{itemize}
\item[$(Q1)$] if $p, q\in {\cal P}$ (possibly $p = q$) then $\langle \varepsilon(p)\cup\varepsilon(q)\rangle\cap
K = 0$. In particular $K\cap \varepsilon(p) = 0$ for every point $p\in{\cal P}$. \item[$(Q2)$] For $p\in{\cal
P}$ and $l\in{\cal L}$, if $p\not\in l$ then $\varepsilon(p)\cap\langle \varepsilon(l)\cup K\rangle = 0$.
\item[$(Q3)$] There exits a nonnegative integer $k < d$ such that
$\mathrm{dim}(K\cap\langle\varepsilon(l)\rangle) = k$ for every line $l\in{\cal L}$.
\end{itemize}

Then the function $\varepsilon/K$ mapping $p\in{\cal P}$ to $[\langle \varepsilon(p)\cup K\rangle/K]$ is a
$(d-k)$-embedding of $\Gamma$ in $\mathrm{PG}(V/K)$ and the canonical projection of $V$ onto $V/K$ is a morphism
from $\varepsilon$ to $\varepsilon/K$. Indeed $(Q1)$ forces $\varepsilon/K$ to be injective, $(Q2)$ implies
$(E2)$ for $\varepsilon/K$ and $(Q3)$ implies $(E1)$ with $d$ replaced by $d-k$ ($> 0$ as $k < d$; note also
that the condition $k < d$ could be removed from $(Q3)$, since it follows from $(Q1)$ and $(Q2)$).

We call $\varepsilon/K$ a $k$-{\em quotient} of $\varepsilon$ (also a {\em quotient} of $\varepsilon$, omitting
to mention $k$), and we say that $K$ {\em defines a $k$-quotient} of $\varepsilon$ (a {\em quotient} of
$\varepsilon$, for short).

If $\varepsilon$ is a projective embedding then $(Q1)$ implies both $(Q2)$ and $(Q3)$. If $\varepsilon$ is a
$1$-embedding then $(Q2)$ does not follow from $(Q1)$ but $(Q1)$ still implies $(Q3)$. Clearly, all quotients of
a $1$-embedding are $0$-quotients.

If $\varepsilon$ is a quadratic embedding defined over a field $\mathbb{F}$ with $\mathrm{char}(\mathbb{F}) \neq
2$, then $(Q1)$ implies $(Q3)$ with $k = 0$. In this case all quotients of $\varepsilon$ are quadratic (see
\cite{CP1}, \cite{CP2}). On the other hand, let $\mathrm{char}(\mathbb{F}) = 2$ and let $K$ be a subspace of $V$
satisfying $(Q1)$ and $(Q2)$. Then $\mathrm{dim}(K\cap\langle\varepsilon(l)\rangle) \leq 1$ for every line
$l\in{\cal L}$ (see \cite{CP1}, \cite{CP2}). Thus, assuming $(Q3)$ is equivalent to assume that $k :=
\mathrm{dim}(K\cap\langle\varepsilon(l)\rangle)$ does not depend on the choice of $l$. If $k = 0$ then
$\varepsilon/K$ is quadratic. If $k = 1$ then $[K]\cap \langle\varepsilon(l)\rangle^{\mathrm{pr}}$ is the
nucleus of the conic $\varepsilon(l)$. In this case $\varepsilon/K$ is a $1$-embedding. It is a projective
embedding if $\mathbb{F}$ is perfect, otherwise it can be lax (see \cite{CP1}, \cite{CP2}).

\subsubsection{Hulls and universality}

A $d$-embedding of $\Gamma$ is said to be {\em relatively universal} (dominant in Tits \cite{Tits} and in
\cite{PasEE}) when it is not a proper morphic image of any other $d$-embedding of $\Gamma$. As we will show in a
few lines, every $d$-embedding $\varepsilon$ of $\Gamma$ admits a {\em hull} $\tilde{\varepsilon}$, uniquely
determined up to isomorphisms by the following properties: $\tilde{\varepsilon}$ is a $d$-embedding of $\Gamma$
and $\tilde{\varepsilon}\geq \varepsilon'$ for every $d$-embedding $\varepsilon'$ of $\Gamma$ such that
$\varepsilon'\geq\varepsilon$. In particular, $\tilde{\varepsilon}\geq \varepsilon$. Clearly,
$\tilde{\varepsilon}$ is relatively universal. A $d$-embedding is relatively universal if and only if it is its
own hull.

The hull $\tilde{\varepsilon}$ of a $d$-embedding $\varepsilon:\Gamma\stackrel{d}{\rightarrow}\mathrm{PG}(V)$
can be constructed as follows, by a word for word rephrasing of the construction given by Ronan \cite{Ronan} for
hulls of projective embeddings. Denoted by $\cal F$ the set of point-line flags of $\Gamma$, consider the
following presheaf:
\[{\bf F}(\varepsilon) = (\{V_x\}_{x\in{\cal P}\cup{\cal L}},\{\iota_{p,l}\}_{(p,l)\in{\cal F}})\]
where if $x\in{\cal P}$ then $V_x = \varepsilon(x)$ (regarded as a 1-dimensional vector space), if $x\in{\cal
L}$ then $V_x = \langle\varepsilon(x)\rangle = \langle V_p\rangle_{p\in x}$ (a $(d+1)$-dimensional vector space)
and $\iota_{p,l}$ is the inclusion embedding of $V_p$ in $V_l$, for every flag $(p,l)\in{\cal F}$.

The vector spaces $V_p$ and $V_l$ for $p\in{\cal P}$ and $l\in{\cal L}$ are taken from $V$, but when forming the
presheaf ${\bf F}(\varepsilon)$ we regard them as abstract $1$- and $(d+1)$-dimensional vector spaces. So, we
can consider their formal direct sum $O(\varepsilon) := \oplus_{x\in{\cal P}\cup{\cal L}}V_x$. Let $J$ be the
subspace of $O(\varepsilon)$ spanned by the vectors $v-\iota_{p,l}(v)$ for every flag $(p,l)\in{\cal F}$ and
every vector $v\in \varepsilon(p)$.

Put $V(\varepsilon) = O(\varepsilon)/J$ and define the mapping $\tilde{\varepsilon}$ from ${\cal P}$ to the set
of $1$-dimensional subspaces of $V(\varepsilon)$ as follows: $\tilde{\varepsilon}(p) = \langle V_p\cup J\rangle$
for every point $p\in{\cal P}$. Then $\tilde{\varepsilon}$ is a $d$-embedding of $\Gamma$ in
$\mathrm{PG}(V(\varepsilon))$ and the natural projection of $V(\varepsilon)$ onto $V$ is a morphism from
$\tilde{\varepsilon}$ to $\varepsilon$.

\begin{proposition}
The embedding $\tilde{\varepsilon}$ defined as above is the hull of $\varepsilon$.
\end{proposition}
The proof is straightforward. We leave it to the reader.

The previous proposition answers a theoretical existence question but it is not of great help when we want to
check if a given embedding is relatively universal or compute the dimension of its hull. When dealing with
projective embeddings questions like these can be answered, sometimes easily, sometimes with some hard work.
They are more difficult for $1$-embeddings in general and become extremely difficult when we turn to quadratic
embeddings. Nevertheless, something can be said in this latter case too. For instance, it is proved in
\cite{TVM} and \cite{TVMmore} that every $2$-embedding of a projective space of finite dimension $n$ with at
least four points on each line is relatively universal (and has dimension equal to ${{n+2}\choose 2}-1$)
provided that its dimension is at least ${{n+2}\choose 2}-1$. The same conclusion, but only for quadratic
embeddings, is obtained in \cite{CP1} by an easier argument (see also \cite{CP2}). One more result in this
trend, taken from \cite{CP2}, will be mentioned in Section \ref{Grassmann-Emb-GQ}, Proposition \ref{W(3,q)
univ}.

So far for hulls and relative universality. Turning to absolute universality, let $\bf C$ be a nonempty class of
$d$-embeddings of $\Gamma$ defined over $\mathbb{F}$, for a given positive integer $d$ and a given division ring
$\mathbb{F}$. We assume that $\bf C$ is closed under isomorphism and under taking $0$-quotients and hulls and
that if $\varepsilon\in{\bf C}$ then $\varepsilon\cdot g\in{\bf C}$ for every automorphism $g$ of $\Gamma$. A
member $\varepsilon$ of $\bf C$ is said to be {\em absolutely universal} (in $\bf C$) if it is the hull of all
members of $\bf C$. Clearly, the absolutely universal member of $\bf C$, if it exists, is uniquely determined up
to isomorphisms.

For the above definition to be interesting the class $\bf C$ must be large enough, but not too large, otherwise
there is no chance for $\bf C$ to admit an absolutely universal member. For instance, $\bf C$ can be chosen as
the class of all projective embeddings of $\Gamma$ defined over a given division ring or the class of all
quadratic embeddings of $\Gamma$ defined over a given field.

Absolutely universal projective embeddings have been thoroughly studied. We refer the reader to Kasikova and
Shult \cite{KS} for a very far-reaching sufficient condition for the existence of the absolutely universal
projective embedding. A bit more on this topic can be found in Blok and Pasini \cite{BP}.

Given a division ring $\mathbb{F}$ and a geometry $\Gamma$ admitting a lax $\mathbb{F}$-$1$-embedding, it can
happen that there is a subring $\mathbb{F}_0$ of $\mathbb{F}$ such that all $\mathbb{F}$-$1$-embeddings of
$\Gamma$ are obtained as scalar extensions of projective embeddings defined over $\mathbb{F}_0$. If this is the
case then we can exploit what we know on absolutely universal projective embeddings of $\Gamma$ to investigate
absolutely universal $1$-embeddings of $\Gamma$.

Almost nothing is known on the existence of absolutely universal quadratic embeddings.

\subsubsection{Homogeneity}\label{Homogeneity}

Given a $d$-embedding $\varepsilon:\Gamma\stackrel{d}{\rightarrow}\mathrm{PG}(V)$ and an automorphism $g$ of
$\Gamma$, a {\em lifting} of $g$ to $\mathrm{PG}(V)$ through $\varepsilon$ is an invertible semilinear mapping
$\hat{g}\in \Gamma\mathrm{L}(V)$ such that $\hat{g}\varepsilon = \varepsilon g$. Clearly, if $\hat{g}$ is a
lifting of $g$ then the coset $\hat{g}\cdot \mathrm{Aut}(\varepsilon)$ is the family of all liftings of $g$. In
particular, if $\mathrm{Aut}(\varepsilon) = Z(\mathrm{GL}(V))$ (as when $\varepsilon$ is projective or
quadratic) then the lifting $\hat{g}$ is unique modulo scalars.

The set of all elements of $\mathrm{Aut}(\Gamma)$ that lift to $\mathrm{PG}(V)$ through $\varepsilon$ is a
subgroup $\mathrm{Aut}_\varepsilon(\Gamma)$ of $\mathrm{Aut}(\Gamma)$. Given a subgroup $G$ of
$\mathrm{Aut}(\Gamma)$, if $G\leq \mathrm{Aut}_\varepsilon(\Gamma)$ we say that $\varepsilon$ is $G$-{\em
homogeneous}. Let $\varepsilon$ be $G$-homogeneous. The {\em lifting} $\widehat{G}$ of $G$ is the subgroup of
$\Gamma\mathrm{L}(V)$ formed by the liftings of the elements of $G$. It contains $\mathrm{Aut}(\varepsilon)$ as
a normal subgroup and $\widehat{G}/\mathrm{Aut}(\varepsilon)\cong G$.

If $\varepsilon$ is $\mathrm{Aut}(\Gamma)$-homogenous then we say it is {\em fully homogeneous}, also just {\em
homogeneous}, for short.

Clearly, the hull of a $G$-homogeneous $d$-embedding is $G$-homogeneous. Conversely, let $K$ define a quotient
of $\varepsilon$. If $\varepsilon$ is $G$-homogeneous and $K$ is stabilized by the lifting of $G$ through
$\varepsilon$ then $\varepsilon/K$ is $G$-homogeneous.

\subsection{Two ways to construct $d$-embeddings}\label{ImmersioniIndotte}

\subsubsection{Veronesean embeddings}\label{veronesean embeddings}

Given a point-line geometry $\Gamma = ({\cal P},{\cal L})$ let $\varepsilon:\Gamma\rightarrow \mathrm{PG}(V)$ be
a projective embedding. If $\eta:\mathrm{PG}(V)\stackrel{d}{\rightarrow}\mathrm{PG}(V')$ is a $d$-embedding of
$\mathrm{PG}(V)$ then the composition $\eta\cdot\varepsilon$ is a $d$-embedding of $\Gamma$ in the subspace of
$\mathrm{PG}(V')$ spanned by $\eta({\cal P})$.

In particular, given a vector space $V$ of finite dimension $n$ over a field $\mathbb{F}$, let $V^\otimes =
(V\otimes V)/\langle x\otimes y - y\otimes x\rangle_{x,y\in V}$ be the symmetrized tensor square of $V$. As
$\mathrm{dim}(V^\otimes) = {{n+1}\choose 2}$, we can choose the pairs $(i,j)$ with $1\leq i\leq j\leq n$ as
indices for the vectors of a basis of $V^\otimes$. Chosen a basis $E = \{e_i\}_{i=1}^n$ of $V$ and a basis
$E^\otimes = \{e_{i,j}\}_{1\leq i\leq j\leq n}$ of $V^\otimes$, the {\em quadric veronesean map} from $V$ to
$V^\otimes$ (relative to $E$ and $E^\otimes$) is the function $\eta^{\mathrm{ver}}$ mapping a vector $x =
\sum_{i=1}^ne_ix_i$ of $V$ onto the vector $\eta^{\mathrm{ver}}(x) = \sum_{i\leq j}e_{i,j}x_ix_j$ of
$V^\otimes$.

The mapping $\eta^{\mathrm{ver}}$ naturally defines a quadratic embedding of $\mathrm{PG}(V)$ in
$\mathrm{PG}(V^\otimes)$, which we still denote by the symbol $\eta^{\mathrm{ver}}$. The composition
$\eta^{\mathrm{ver}}\varepsilon$ is a quadratic embedding of $\Gamma$ in a subspace of $\mathrm{PG}(V^\otimes)$.
We call it the {\em veronesean embedding} of $\Gamma$ {\em induced} by $\varepsilon$ and we denote it by the
symbol $\varepsilon^{\mathrm{ver}}$.

The image ${\cal V} = \eta^{\mathrm{ver}}(\mathrm{PG}(V))$ of the set of points of $\mathrm{PG}(V)$ by
$\eta^{\mathrm{ver}}$ is a well known projective variety of dimension $n-1$, called the {\em veronesean
variety}, described by the following set of equations in the unknowns $x_{i,j}$ for $1\leq i \leq j \leq n$,
where we put $x_{j,i} := x_{i,j}$ when $j > i$:

\begin{equation}\label{Veronese-equazioni}
\left.\begin{array}{cl}
x_{i,i}x_{j,j} = x_{i,j}^2 & \mbox{for $i < j$},\\
x_{i,j}x_{i,k} = x_{i,i}x_{j,k} & \mbox{for $j < k$ and $i \neq j, k$},\\
x_{i,j}x_{k,h} = x_{i,k}x_{j,h} & \mbox{for $i < j < h, k$ and $h\neq k$}.
\end{array}\right\}
\end{equation}
(See Hirschfeld and Thas \cite{HT}.) Clearly ${\cal V}\supseteq \varepsilon^{\mathrm{ver}}({\cal P})$. In many
interesting cases $\varepsilon^{\mathrm{ver}}({\cal P})$ is a subvariety of $\cal V$. For instance, if
$\varepsilon({\cal P})$ is a quadric in $\mathrm{PG}(V)$ then $\varepsilon^{\mathrm{ver}}({\cal P})$ is a
hyperplane section of $\cal V$.

The embedding $\eta^{\mathrm{ver}}$ is homogeneous. Therefore, if $\varepsilon$ is $G$-homogeneous for a
subgroup $G$ of $\mathrm{Aut}(\Gamma)$ then $\varepsilon^{\mathrm{ver}}$ is $G$-homogeneous.

More generally, let $\varepsilon_1:\Gamma\rightarrow\mathrm{PG}(V_1)$ and
$\varepsilon_2:\Gamma\rightarrow\mathrm{PG}(V_2)$ be two projective embeddings of $\Gamma$ and $f:V_1\rightarrow
V_2$ a morphism from $\varepsilon_1$ to $\varepsilon_2$. A unique semilinear mapping
$f^\otimes:V_1^{\otimes}\rightarrow V^\otimes_2$ exists such that $f^\otimes \eta_1^{\mathrm{ver}} =
\eta_2^{\mathrm{ver}} f$, where $\eta_1^{\mathrm{ver}}$ and $\eta_2^{\mathrm{ver}}$ are the quadric veronesean
mappings of $V_1$ in $V_1^\otimes$ and $V_2$ in $V_2^\otimes$ respectively. Clearly, $f^\otimes$ is a morphism
from $\varepsilon_1^{\mathrm{ver}}$ to $\varepsilon_2^{\mathrm{ver}}$.

\subsubsection{Grasmmann embeddings}\label{Grassmann embeddings}

Given a point line geometry $\Gamma = ({\cal P},{\cal L})$, for every point $p\in{\cal L}$ let $L(p)$ be the set
of lines of $\Gamma$ through $p$. Put ${\cal P}^* = {\cal L}$, ${\cal L}^* = \{L(p)\}_{p\in{\cal P}}$ and
$\Gamma^* = ({\cal P}^*,{\cal L}^*)$. Then $\Gamma^*$ is a point-line geometry, isomorphic to the dual $({\cal
L},{\cal P})$ of $\Gamma$.

Given a field $\mathbb{F}$ and a vector space $V$ over $\mathbb{F}$ of finite dimension at least $4$, suppose
that a projective embedding $\varepsilon:\Gamma^*\rightarrow \mathrm{PG}(V)$ is given such that

\begin{itemize}
\item[$(E^*_d)$] $\mathrm{dim}(\langle\varepsilon(L(p))\rangle^{pr}_{p\in l}) = d+1$ for every point $l$ of
$\Gamma^*$ (line of $\Gamma$).
\end{itemize}

We recall that the {\em line-grassmannian} $\mathrm{Gr}_2(\mathrm{PG}(V))$ of $\mathrm{PG}(V)$ is the point-line
geometry where the lines of $\mathrm{PG}(V)$ are taken as points and the lines are the sets of lines of
$\mathrm{PG}(V)$ incident to a point-plane flag of $\mathrm{PG}(V)$. Put $V^\wedge = V\wedge V$ and let
$\iota_2$ be the function mapping every line $\langle [x],[y]\rangle^{\mathrm{pr}}$ of $\mathrm{PG}(V)$ to the
point $[x\wedge y]$ of $\mathrm{PG}(V^\wedge)$. It is well known that $\iota_2$ is a projective embedding of
$\mathrm{Gr}_2(\mathrm{PG}(V))$ in $\mathrm{PG}(V^\wedge)$. We call it the {\em natural embedding} of
$\mathrm{Gr}_2(\mathrm{PG}(V))$.

For every point $p\in{\cal P}$ of $\Gamma$ put $\varepsilon^{\mathrm{gr}}(p) = \iota_2(\varepsilon(L(p)))$.

\begin{proposition}
The mapping $\varepsilon^{\mathrm{gr}}$ defined as above is a $d$-embedding of $\Gamma$ in the subspace of
$\mathrm{PG}(V^\wedge)$ spanned by $\iota_2({\cal L}^*)$ ($= \iota_2({\cal P})$).
\end{proposition}
{\bf Proof.} The mapping $\varepsilon^{\mathrm{gr}}$ is injective, since both $\varepsilon$ and $\iota_2$ are
injective and distinct points of $\Gamma$ are incident with distinct sets of lines of $\Gamma$. Property $(E3)$
trivially holds, since we have chosen $\langle \iota_2({\cal L}^*)\rangle^{\mathrm{pr}}$ as the codomain of
$\varepsilon^{\mathrm{gr}}$.

Every line $l\in{\cal L}$ of $\Gamma$ is mapped by $\varepsilon^{\mathrm{gr}}$ onto the set
$\{\iota_2(\varepsilon(L(p)))\}_{p\in l}$. By $(E^*_d)$ the set $X = \{\varepsilon(L(p))\}_{p\in l}$ spans a
$d+1$-dimensional subspace $\langle X\rangle^{\mathrm{pr}}$ of $\mathrm{PG}(V)$ containing the point
$\varepsilon(l)$. On the other hand, it is well known that, for every point $a$ of $\mathrm{PG}(V)$, denoted by
$\mathrm{St}(a)$ the set of lines of $\mathrm{PG}(V)$ through $a$, the set $\iota_2(\mathrm{St}(a))$ is a
subspace of $\mathrm{PG}(V^\wedge)$ and $\iota_2$ induces an isomorphism from the residue of $a$ in
$\mathrm{PG}(V)$ to the subspace $\iota_2(\mathrm{St}(a))$. It follows that $\langle
\varepsilon^{\mathrm{gr}}(l)\rangle^{\mathrm{pr}} = \langle
\iota_2(\varepsilon(L(p)))\rangle^{\mathrm{pr}}_{p\in l}$ is a $d$-dimensional subspace of
$\mathrm{PG}(V^\wedge)$. So, $\varepsilon^{\mathrm{gr}}$ satisfies $(E1)$.

Finally, $\langle \varepsilon^{\mathrm{gr}}(l)\rangle^{\mathrm{pr}}\cap \iota_2(\mathrm{PG}(V)) = \iota_2(\mathrm{St}(\varepsilon(l)))$. It follows that $\langle \varepsilon^{\mathrm{gr}}(l)\rangle^{\mathrm{pr}}\cap \varepsilon^{\mathrm{gr}}({\cal P}) = \varepsilon^{\mathrm{gr}}(l)$ (recall that $\varepsilon$ maps lines of $\Gamma^*$ surjectively onto lines of $\mathrm{PG}(V)$, since it is a projective embedding). So, $\varepsilon^{\mathrm{gr}}$ also satisfies $(E2)$.  \hfill $\Box$

\bigskip

We call $\varepsilon^{\mathrm{gr}}$ the {\em Grassmann embedding} of $\Gamma$ {\em induced} by $\varepsilon$,
also a {\em Grassmann embedding} for short.

The $\iota_2$-image ${\cal G} := \iota_2(\mathrm{Gr}_2(\mathrm{PG}(V)))$ of the set of points of
$\mathrm{Gr}_2(\mathrm{PG}(V))$ (lines of $\mathrm{PG}(V)$) is a well known projective variety, called the {\em
line-Grassmann variety} of $\mathrm{PG}(V)$.
 Clearly, $\varepsilon^{\mathrm{gr}}({\cal P})$ is a subset of ${\cal G}$.
In many interesting cases it is a subvariety of $\cal G$.

The embedding $\iota_2$ is homogeneous. Moreover $\mathrm{Aut}(\Gamma)$ and $\mathrm{Aut}(\Gamma^*)$ are
canonically isomorphic. Therefore, if $\varepsilon$ is $G^*$-homogeneous for a subgroup $G^*$ of
$\mathrm{Aut}(\Gamma^*)$ and $G$ is the subgroup of $\mathrm{Aut}(\Gamma)$ corresponding to $G^*$ then
$\varepsilon^{\mathrm{gr}}$ is $G$-homogeneous.

More generally, let $\varepsilon_1:\Gamma^*\rightarrow\mathrm{PG}(V_1)$ and
$\varepsilon_2:\Gamma^*\rightarrow\mathrm{PG}(V_2)$ be two projective embeddings of $\Gamma^*$ satisfying
condition $(E^*_d)$, for possibly different choices of $d$. Let $f:V_1\rightarrow V_2$ be a morphism from
$\varepsilon_1$ to $\varepsilon_2$. There exists a semilinear mapping $f^\wedge:V_1^\wedge\rightarrow
V_2^\wedge$, unique modulo scalars, which commutes with the natural embeddings of
$\mathrm{Gr}_2(\mathrm{PG}(V_1))$ and $\mathrm{Gr}(\mathrm{PG}(V_2))$. Clearly, $f^\wedge$ is morphism from
$\varepsilon_1^{\mathrm{ver}}$ to $\varepsilon_2^{\mathrm{ver}}$.

\subsubsection{A more general setting}

The previous definition of $\varepsilon^{\mathrm{gr}}$ is a special case of a more general setting, quite
familiar to everyone working on embeddings of dual structures, as dual polar spaces for instance. We recall it
below, even if we will make no use of it in this paper.

Let $\Gamma_0$ be a point-line geometry equipped with two families of subspaces ${\cal P}^*$ and ${\cal L}^*$
such that every member of ${\cal P}^*$ is properly contained in some members of ${\cal L}^*$ and in no member of
${\cal P}^*$, every member of ${\cal L}^*$ properly contains some members of ${\cal P}^*$ but no member of
${\cal L}^*$ and the pair $\Gamma^* := ({\cal P}^*,{\cal L}^*)$ equipped with inclusion as the incidence
relation is a rank $2$ geometry isomorphic to the dual $({\cal L},{\cal P})$ of $\Gamma$. Chosen an isomorphism
from the dual of $\Gamma$ to $\Gamma^*$ let $\lambda$ be the bijection from ${\cal P}$ to ${\cal L}^*$ induced
by that isomorphism. Let $\varepsilon:\Gamma_0\rightarrow\mathrm{PG}(V)$ be a projective embedding such that for
a given integer $k\geq 2$ and any choice of $X\in{\cal P}^*$ and $Y\in{\cal L}^*$ we have $\mathrm{dim}(\langle
\varepsilon(X)\rangle^{\mathrm{pr}}) = k-2$, $\mathrm{dim}(\langle\varepsilon(Y)\rangle^{\mathrm{pr}}) = k-1$
and $\langle\varepsilon(X)\rangle^{\mathrm{pr}}\subseteq\langle\varepsilon(Y)\rangle^{\mathrm{pr}}$ if and only
if $X\subseteq Y$. Then, if $\iota_k$ is the usual projective embedding of the $k$-grassmannian of
$\mathrm{PG}(V)$ in $\mathrm{PG}(\wedge^kV)$ mapping a $(k-1)$-subspace $\langle [x_1],...,
[x_k]\rangle^{\mathrm{pr}}$ of $\mathrm{PG}(V)$ onto the point $[\wedge_{i=1}^k x_i]$, the composition
$\varepsilon^{\mathrm{gr}} := \iota_k\varepsilon\lambda$, which maps $p\in{\cal P}$ onto
$\iota_k(\varepsilon(\lambda(p)))$, is an injective mapping from the point-set $\cal P$ of $\Gamma$ to the set
of points of $\mathrm{PG}(\wedge^kV)$.

Suppose moreover that for every $X\in{\cal P}^*$ the set $\{\varepsilon(Y)\}_{X\subset Y\in{\cal L}^*}$ spans a
$(k-1+d)$-dimensional subspace of $\mathrm{PG}(V)$. Then $\varepsilon^{\mathrm{gr}}$ is a $d$-embedding of
$\Gamma$ in a subspace of $\mathrm{PG}(\wedge^kV)$.

\section{Line grassmannians of polar spaces}\label{Section 3}

Throughout this section $V$ is a vector space of finite dimension $n \geq 4$ over a field $\mathbb{F}$,
$V^\wedge := V\wedge V$ and $\iota_2$ is the natural embedding of $\mathrm{Gr}_2(\mathrm{PG}(V))$ in
$\mathrm{PG}(V^\wedge)$, as in Subsection \ref{Grassmann embeddings}. As in that subsection, ${\cal G} :=
\iota^{\mathrm{gr}}_2(\mathrm{Gr}_2(\mathrm{PG}(V)))$, the line-Grassmann variety in $\mathrm{PG}(V^\wedge)$.

Chosen an ordered basis $E = (e_1,..., e_n)$ of $V$, the vectors $e_{i,j} := e_i\wedge e_j$ for $1\leq i < j
\leq n$ form a basis $E^\wedge$ of $V^\wedge$. We call $E^\wedge$ the basis of $V^\wedge$ {\em canonically
associated} to $E$. The embedding $\iota_2$ maps the line $\langle [\sum_{i=1}^ne_ix_i],
[\sum_{i=1}^ne_iy_i]\rangle^{\mathrm{pr}}$ of $\mathrm{PG}(V)$ onto the point $[\sum_{i<j}e_{i,j}x_{i,j}]$ of
$\mathrm{PG}(V^\wedge)$, where $x_{i,j} = x_iy_j-x_jy_i$ for any choice of $i < j$.

\subsection{Properties of $\cal G$}\label{properties grassmannian}

In this subsection we recall a few well known properties of $\cal G$. This variety is described by the following
set of equations (Hirschfeld and Thas \cite{HT}):

\begin{equation}\label{Grassmann-equazioni}
x_{i,j}x_{k,h} - x_{i,k}x_{j,h}+ x_{i,h}x_{j,k} = 0, \hspace{5 mm} (1\leq i < j < k < h\leq n).
\end{equation}
Given any two non-proportional vectors $a = \sum_{i=1}^ne_ia_i$ and $b = \sum_{i=1}^ne_ib_i$ of $V$ let $c =
\sum_{i<j}e_{i,j}c_{i,j} = a\wedge b$, namely $c_{i,j} = a_ib_j-a_jb_i$ for $1\leq i < j\leq n$. The tangent
space $\mathrm{Tan}({\cal G})_{[c]}$ of $\cal G$ at the point $[c]$ is the subspace of $V^\wedge$ described by
the following set of linear equations, for $1\leq i < j < k < h\leq n$:
\begin{equation}\label{Grassmann-equazioni-tangente}
c_{k,h}x_{i,j} + c_{i,j}x_{k,h} - c_{j,h}x_{i,k} - c_{i,k}x_{j,h} + c_{j,k}x_{i,h} + c_{i,h}x_{j,k} = 0.
\end{equation}
This linear system contains $n\choose 4$ equations but it has rank equal to ${{n-2}\choose 2}$, independently of
the choice of the point $[c]\in {\cal G}$ (Hirschfeld and Thas \cite{HT}). Accordingly,
$\mathrm{dim}(\mathrm{Tan}({\cal G}_{[c]})) = 2n-3 = {n\choose 2}-{{n-2}\choose 2}$. Hence $\cal G$ has
dimension $\mathrm{dim}({\cal G}) = 2n-4$.

The tangent space $\mathrm{Tan}({\cal G})_{[c]}$ also admits another description, as we shall see in a few lines.

Given a non-zero vector $a = \sum_{i=1}^ne_ia_i$ of $V$, let $\mathrm{St}([a])$ be the set of lines of
$\mathrm{PG}(V)$ containing the point $[a]$, as in Subsection \ref{Grassmann embeddings}. Then $S_{[a]} :=
\iota_2(\mathrm{St}([a]))$ is a subspace of $\mathrm{PG}(V^\wedge)$ and it is described by the following
$n\choose 3$ linear equations:
\begin{equation}\label{Grassmann-equazioni-star}
a_ix_{j,k} - a_jx_{i,k} + a_kx_{i,j} = 0, \hspace{5 mm} (1\leq i < j < k\leq n).
\end{equation}
This linear system has rank equal to ${n-1}\choose 2$. Indeed, suppose that $a_1\neq 0$, to fix ideas. Then the
equations (\ref{Grassmann-equazioni-star}) with $i = 1$ form a maximal independent subset of the whole set
(\ref{Grassmann-equazioni-star}). Accordingly $S_{[a]}$, regarded as a subspace of $V^\wedge$, has vector
dimension $\mathrm{dim}(S_{[a]}) = n-1$.

Given another point $[b]\neq [a]$ in $\mathrm{PG}(V)$, let $c = a\wedge b$. Then $S_{[a]}\cap S_{[b]} = [c]$.
Hence, regarded $S_{[a]}$ and $S_{[b]}$ as subspaces of $V^\wedge$, the sum $S_{[a]}+ S_{[b]}$ has dimension
equal to $2(n-1)-1 = 2n-3 = \mathrm{dim}(\mathrm{Tan}({\cal G})_{[c]})$. It is not difficult to check that both
$S_{[a]}$ and $S_{[a]}$ are contained in $\mathrm{Tan}({\cal G})_{[c]}$. Therefore $\mathrm{Tan}({\cal G})_{[c]}
= S_{[a]}+S_{[b]}$.

\subsection{A matrix notation}

For every vector $v = \sum_{i<j}e_{i,j}v_{i,j}$ of $V^\wedge$ let
\begin{equation}\label{Av}
A_v = (v_{i,j})_{i,j=1}^n
\end{equation}
where we put $v_{i,i} = 0$ for every $i$ and $v_{i,j} = - v_{j,i}$ for $i > j$. So, $A_v$ is an anti-symmetric
matrix. The function $\alpha$ mapping $v\in V^\wedge$ to $\alpha(v) := A_v$ is an isomorphism from $V^\wedge$ to
the space of anti-symmetric $n\times n$ matrices over $\mathbb{F}$. In this way the points of ${\cal G}$ are
represented by matrices of the form $xy^T-yx^T$ for $x, y\in V$ with $\mathrm{dim}(\langle x,y\rangle) = 2$.
(Needless to say, when writing $xy^T-yx^T$ we regard the vectors of $V$ as $(n\times 1)$-matrices with the
$E$-coordinates as the entries.)

We can dualize the above as follows. The linear functionals of $V$ form the dual $V^*$ of $V$. Let $E^* =
(e^*_1,..., e^*_n)$ be the basis of $V^*$ dual of $E$, namely $e^*_i(e_j) = \delta_{i,j}$ (Kronecker symbol) for
every choice of $i,j = 1, 2,..., n$. The vectors $e^*_{i,j} := e^*_i\wedge e^*_j$ for $i < j$ form the basis
$E^{*\wedge}$ of $V^{*\wedge} := V^*\wedge V^*$ canonically associated to $E^*$. Given a vector $\xi =
\sum_{i<j}e^*_{i,j}\xi_{i,j}$ of $V^{*\wedge}$ we put $A^*_\xi = (\xi_{i,j})_{i,j=1}^n$, with $\xi_{i,i} = 0$
and $\xi_{i,j} = -\xi_{j,i}$ if $i > j$, as in (\ref{Av}). Thus $A^*_\xi$ is an anti-symmetric matrix and the
clause $\alpha^*(\xi) = A^*_\xi$ defines an isomorphism $\alpha^*$ from $V^{*\wedge}$ to the vector space of
anti-symmetric $n\times n$ matrices. In particular, if $\theta$ and $\zeta$ are non-proportional vectors of
$V^*$ then $\alpha^*(\theta\wedge \zeta) = \theta\zeta^T-\zeta\theta^T$, where $\theta$ and $\zeta$ are regarded
as $n\times 1$-matrices with the $E^*$-coordinates as the entries.

The equations gathered in (\ref{Grassmann-equazioni}) can be rephrased as properties of anti-symmetric matrices,
as follows: a vector $v\in V^\wedge$ satisfies (\ref{Grassmann-equazioni}) for all choices of $i < j < k < h$ if
and only if all principal $(4\times 4)$-submatrices of $A_v$ are singular. Admittedly, this rephrasing is not a
great improvement, but anti-symmetric matrices are useful in other contexts, as the description of the lifting
of a linear mapping of $V$ to $V^\wedge$ and the characterization of the inclusion of a line of $\mathrm{PG}(V)$
in a subspace of codimension $2$.

\subsubsection{Lifting linear mappings from $V$ to $V^\wedge$}

The embedding $\iota_2$ is homogeneous. In particular, let $M$ be the matrix representing $f\in \mathrm{GL}(V)$
with respect to $E$ by a matrix $M$, namely $f(x) = Mx$ for every $x\in V$. Then
\[\begin{array}{ccccc}
A_{f(x)\wedge f(y)} & = & Mx(My)^T- My(Mx)^T & = & \\
 & = & Mxy^TM^T-Myx^TM^T & = & MA_{x\wedge y}M^T
\end{array}\]
for any two non-proportional vectors $x, y\in V$. So, if $\mathrm{PG}(f)$ is the automorphism of
$\mathrm{PG}(V)$ represented by $f$, then the linear mapping $f^\wedge$ of $V^\wedge$ defined by the clause
\begin{equation}\label{action}
f^\wedge(v) = \alpha^{-1}(MA_vM^T) ~~~ \mbox{(for all $v\in V^\wedge$)}
\end{equation}
is a representative of the lifting of $\mathrm{PG}(f)$ to $\mathrm{PG}(V^\wedge)$. We call $f^\wedge$ the {\em
lifting} of $f$ to $V^\wedge$.

The mapping $f$ induces on $V^*$ the linear mapping $f^*$ acting as follows: $f^*(\theta) = M^{-T}\theta$ for
every $\theta\in V^*$. We have
\[\begin{array}{rl}
A^*_{f^*(\theta)\wedge f^*(\zeta)} & =  M^{-T}\theta(M^{-T}\zeta)^T- M^{-T}\zeta(M^{-T}\theta)^T  =  \\
 & =  M^{-T}\theta\zeta^TM^{-1}-M^{-T}\zeta\theta^TM^{-1} ~~ = ~~ M^{-T}A^*_{\theta\wedge \zeta}M^{-1}
\end{array}\]
for any two non-proportional vectors $\theta, \zeta\in V^*$. Thus, the linear mapping $f^{*\wedge}$ of
$V^{*\wedge}$ defined by the clause
\begin{equation}\label{action-duale}
f^{*\wedge}(\xi) = {\alpha^*}^{-1}(M^{-T}A^*_\xi M^{-1})~~~\mbox{(for all $\xi\in V^{*\wedge}$)}
\end{equation}
is the lifting of $f^*$ to $V^{*\wedge}$.

\subsubsection{Inclusion of a line in a dual line}

The subspaces of $\mathrm{PG}(V)$ of codimension $2$, also called {\em dual lines} of $\mathrm{PG}(V)$,
naturally correspond to the lines of $\mathrm{PG}(V^*)$. Explicitly, recalling that the vectors of $V^*$ are the
linear functionals of $V$, a line $\langle [\theta], [\zeta]\rangle^{\mathrm{pr}}$ of $\mathrm{PG}(V^*)$
corresponds to the subspace $[\mathrm{Ker}(\theta)\cap\mathrm{Ker}(\zeta)]$ of $\mathrm{PG}(V)$.

\begin{lemma}\label{fundamental lemma}
Let $x, y\in V$ and $\theta, \zeta\in V^*$ be such that $\mathrm{dim}(\langle x,y\rangle) = \mathrm{dim}(\langle
\theta, \zeta\rangle) = 2$. Then we have $\langle x,y\rangle \subseteq
\mathrm{Ker}(\theta)\cap\mathrm{Ker}(\zeta)$ if and only if
\begin{equation}\label{fundamental}
A^*_{\theta\wedge \zeta}A_{x\wedge y} = O
\end{equation}
where $O$ stands for the null square matrix of order $n$.
\end{lemma}
{\bf Proof.} Put $X = A_{x\wedge y}$ and $\Theta = A^*_{\theta\wedge \zeta}$, for short. Then
\[\begin{array}{rcl}
\Theta X  & = & (\theta\zeta^T-\zeta\theta^T)(xy^T-yx^T) = \\
 & = & \theta(\zeta^Tx)y^T-\zeta(\theta^Tx)y^T-\theta(\zeta^Ty)x^T+\zeta(\theta^Ty)x^T = \\
 & = & (\zeta^Tx)\theta y^T-(\theta^Tx)\zeta y^T-(\zeta^Ty)\theta x^T+(\theta^Ty)\zeta x^T.
\end{array}\]
(Recall that $\zeta^Tx$, $\theta^Tx$, $\zeta^Ty$ and $\theta^Ty$ are scalars.) Suppose that $x, y \in
\mathrm{Ker}(\theta)\cap\mathrm{Ker}(\zeta)$. In matrix notation, $\theta^Tx = \zeta^Tx = \theta^Ty = \zeta^Ty =
0$. However $\Theta X = (\zeta^Tx)\theta y^T-(\theta^Tx)\zeta y^T-(\zeta^Ty)\theta x^T+(\theta^Ty)\zeta x^T$, as
shown above. Hence $\Theta X = O$.

Conversely, let $\Theta X = O$. Recall that, if $x = \sum_{i=1}^ne_ix_i$, $y = \sum_{i=1}^ne_iy_i$, $\theta =
\sum_{i=1}^ne^*_i\theta_i$ and $\zeta = \sum_{i=1}^ne^*_i\zeta_i$, then $x_{i,j} = x_iy_j-x_jy_i$ is the
$(i,j)$-entry of $X$ and $\zeta_{i,j} = \theta_i\zeta_j-\theta_j\zeta_i$ is the $(i,j)$-entry of $\Theta$.

By (\ref{action}) and (\ref{action-duale}), if we replace $x$ and $y$ by $f(x)$ and $f(y)$ for a linear mapping
$f\in \mathrm{GL}(V)$ then $\Theta$ is replaced by $\Theta' := M^{-T}\Theta M^{-1}$ and $X$ by $X' := MXM^T$,
where $M$ is the representative matrix of $f$. Accordingly, the product $\Theta X$ is replaced by $\Theta'X' =
M^{-T}\Theta XM^T$. As $M$ is non-singular, we have $\Theta X = O$ if and only if $\Theta'X' = O$. Moreover,
$\langle x,y\rangle\subseteq \mathrm{Ker}(\theta)\cap\mathrm{Ker}(\zeta)$ if and only if $\langle f(x),
f(y)\rangle\subseteq \mathrm{Ker}(f^*(\theta))\cap\mathrm{Ker}(f^*(\zeta))$. We can always choose $f$ in such a
way that $f(x) = e_1$ and $f(y) = e_2$.

By the above, we can assume without loss that $x = e_1$ and $y = e_2$. With this choice of $x$ and $y$ we have
$x_{1,2} = -x_{2,1} = 1$ and $x_{i,j} = 0$ for $\{i,j\}\neq \{1,2\}$. Thus the equation $\Theta X = O$ forces
$\zeta_{i,1} = \zeta_{i,2} = 0$ for every $i = 1, 2,..., n$, namely $\theta_i\zeta_1-\theta_1\zeta_i =
\theta_i\zeta_2-\theta_2\zeta_i = 0$ for every $i =1, 2,..., n$. Suppose that at least one of $\theta_1,
\zeta_1, \theta_2$ or $\zeta_2$ is non-zero, say $\theta_1\neq 0$. Then $\zeta_i =
\frac{\zeta_1}{\theta_1}\theta_i$ for every $i$, namely $\zeta = \frac{\zeta_1}{\theta_1}\theta$, contrary to
the hypotheses of the lemma. It follows that $\theta_1 = \theta_2 = \zeta_1 = \zeta_2 = 0$, namely $\theta^Te_1
= \theta^Te_2 = \zeta^Te_1 = \zeta^Te_2 = 0$. So, $e_1,e_2\in \mathrm{Ker}(\theta)\cap\mathrm{Ker}(\zeta)$,
namely $x, y\in \mathrm{Ker}(\theta)\cap\mathrm{Ker}(\zeta)$. \hfill $\Box$

\subsection{Lines totally isotropic for a sesquilinear form}\label{form-section}

Let $\sigma\in \mathrm{Aut}(\mathbb{F})$ and $\varepsilon\in\mathbb{F}\setminus\{0\}$ be such that $\sigma^2 =
\mathrm{id}_{\mathbb{F}}$ and $\varepsilon^\sigma\varepsilon = 1$. Let $\varphi:V\times V\rightarrow\mathbb{F}$
be a reflexive $(\sigma,\varepsilon)$-sequilinear form on $V$, represented as follows in matrix notation with
respect to the basis $E$ of $V$:
\begin{equation}\label{form}
\varphi(x,y) = (x^{\sigma})^T\Phi y
\end{equation}
for an $n\times n$ matrix $\Phi = (\varphi_{i,j})_{i,j=1}^n$ such that $\Phi^T = \varepsilon \Phi^\sigma$, where
$\Phi^\sigma := (\varphi^\sigma_{i,j})_{i,j=1}^n$. We can rewrite (\ref{form}) as follows. Consider the
semi-linear mapping $\lambda_\varphi:V\rightarrow V^*$ mapping $x$ to the linear functional $\lambda_\varphi(x)$
defined by the clause $\lambda_\varphi(x)(y) = (x^{\sigma})^T\Phi y$ for every $y\in V$, namely
$\lambda_\varphi(x) = \Phi^Tx^\sigma = \varepsilon \Phi^\sigma x^\sigma$ with respect to the basis $E^*$ of
$V^*$. Then (\ref{form}) is equivalent to $\lambda_\varphi(x)(y) = 0$.

We recall that, denoted by $\perp$ the orthogonality relation with respect to $\varphi$, the {\em radical}
$\mathrm{Rad}(\varphi)$ of $\varphi$ is the subspace $V^{\perp}$ of $V$. The form $\varphi$ is said to be {\em
degenerate} if $\mathrm{Rad}(\varphi) \neq 0$. Clearly, $\mathrm{Rad}(\varphi) = \mathrm{Ker}(\lambda_\varphi) =
\mathrm{Ker}(\Phi^\sigma) = (\mathrm{Ker}(\Phi))^\sigma$. So, $\lambda_\varphi$ is an isomorphism if and only if
$\varphi$ is non-degenerate if and only if $\Phi$ is non-singular.

Let $S = \langle x, y\rangle$ be a $2$-dimensional subspace of $\mathrm{PG}(V)$. Then $\lambda_\varphi(S)$ is a
subspace of $V^*$ of dimension at most $2$. Suppose firstly that $\mathrm{dim}(\lambda_\varphi(S)) = 2$ (which
is always the case when $\varphi$ is non-degenerate). Then
\[\begin{array}{rcl}
A^*_{\lambda_\varphi(x)\wedge \lambda_\varphi(y)} & = & \varepsilon^2(\Phi^\sigma x^\sigma(\Phi^\sigma y^\sigma)^T - \Phi^\sigma y^\sigma(\Phi^\sigma x^\sigma)^T) = \\
& = & \varepsilon(\Phi^\sigma (xy^T)^\sigma\Phi - \Phi^\sigma (yx^T)^\sigma\Phi) = ~~ \varepsilon\Phi^\sigma A_{x\wedge y}^\sigma \Phi.\end{array}\]
Namely,
\begin{equation}\label{form-wedge}
A^*_{\lambda_\varphi(x)\wedge \lambda_\varphi(y)}  = \varepsilon\Phi^\sigma A_{x\wedge y}^\sigma\Phi.
\end{equation}
Suppose now that $\mathrm{dim}(\lambda_\varphi(S)) \leq 1$. In this case $\varphi$ is degenerate and $S$ either
intersects the radical $\mathrm{Rad}(\varphi)$ of $\varphi$ in a $1$-dimensional subspace or is entirely
contained in $\mathrm{Rad}(\varphi)$. In either case $A^*_{\lambda_\varphi(x)\wedge \lambda_\varphi(y)} = O$. On
the other hand $\mathrm{Rad}(\varphi) = \mathrm{Ker}(\Phi)$. Hence either $\Phi x = 0$ or $\Phi y = t\Phi x$ for
a scalar $t$ (possibly $t = 0$). In each of these two cases $\Phi xy^T\Phi^T - \Phi yx^T\Phi^T = O$, namely
$\Phi A_{x\wedge y}\Phi^\sigma = O$ (recall that $\Phi^T = \varepsilon \Phi^\sigma$). However $\Phi A_{x\wedge
y}\Phi^\sigma = (\Phi^\sigma A_{x\wedge y}^\sigma \Phi)^\sigma$ (recall that $\sigma^2 =
\mathrm{id}_{\mathbb{F}}$). So, equation (\ref{form-wedge}) remains valid in this case too.

\begin{theorem}\label{main-theorem}
Let $S = \langle x,y\rangle$ and $S' = \langle u,v\rangle$ be $2$-dimensional linear subspaces of $V$. Suppose
that $S\cap \mathrm{Rad}(\varphi) = 0$. Put $X = A_{x\wedge y}$ and $Y = A_{u\wedge v}$. Then $S' \subseteq
S^{\perp}$ if and only if
\begin{equation}\label{form-equation-1}
\Phi^\sigma X^\sigma\Phi Y = O.
\end{equation}
In particular, $S$ is totally isotropic (namely $S\subseteq S^{\perp}$) if and only if
\begin{equation}\label{form-equation-2}
\Phi^\sigma X^\sigma\Phi X = O.
\end{equation}
\end{theorem}
{\em Proof.} The main claim of the theorem immediately follows from equations (\ref{fundamental}) and
(\ref{form-wedge}) while the second claim is just a special case of the first one. The hypothesis
$S\cap\mathrm{Rad}(\varphi) = 0$ is needed in order to apply (\ref{fundamental}). Indeed in the present context
$\lambda_\varphi(x)$ and $\lambda_\varphi(y)$ play the role of $\theta$ and $\zeta$ of Lemma \ref{fundamental
lemma}, but in that lemma we assume that $\mathrm{dim}(\langle \theta, \zeta\rangle) = 2$.  \hfill $\Box$

\begin{corollary}\label{main-corollary}
Suppose that $\varphi$ is non-degenerate. Let $S = \langle x,y\rangle$ and $S' = \langle u,v\rangle$ be
$2$-dimensional linear subspaces of $V$. Put $X = A_{x\wedge y}$ and $Y = A_{u\wedge v}$, as in Theorem
\ref{main-theorem}. Then $S' \subseteq S^{\perp}$ if and only if
\begin{equation}\label{form-equation-1bis}
X^\sigma\Phi Y = O.
\end{equation}
In particular, $S$ is totally isotropic if and only if
\begin{equation}\label{form-equation-2bis}
X^\sigma\Phi X = O.
\end{equation}
\end{corollary}
{\bf Proof.} Trivial from Theorem \ref{main-theorem}, recalling that $\Phi$ is non-singular since $\varphi$ is
non-degenerate. \hfill $\Box$

\begin{corollary}\label{main-corollary-bis}
Let $S = \langle x,y\rangle$ be a $2$-dimensional subspace of $V$ and put $X = A_{x\wedge y}$. We have $S\cap
\mathrm{Rad}(\varphi) \neq 0$ if and only if
\begin{equation}\label{form-wedge-radical}
\Phi X \Phi^\sigma = O.
\end{equation}
\end{corollary}
{\bf Proof.} The `only if' part of this claim is implicit in the comments after formula (\ref{form-wedge}).
Conversely, suppose that (\ref{form-wedge-radical}) holds. Then $\Phi^\sigma X^\sigma \Phi = O$, hence
$\Phi^\sigma X^\sigma \Phi Y = O$ for every $n\times n$ matrix $Y$. By way of contradiction, suppose that $S\cap
\mathrm{Rad}(\varphi) = 0$. Since $\Phi^\sigma X^\sigma \Phi Y = O$ for any $Y$, Theorem \ref{main-theorem}
implies that $S'\subseteq S^\perp$ for every $2$-space $S'$ of $V$. But this is impossible. Indeed
$\mathrm{dim}(S^\perp) = n-2$ since $\mathrm{dim}(S) = 2$ and $S\cap\mathrm{Rad}(\varphi) = 0$. Therefore
$S\cap\mathrm{Rad}(\varphi) \neq 0$. \hfill $\Box$

\begin{corollary}\label{main-corollary-ter}
Let $S = \langle x,y\rangle$ be a $2$-dimensional subspace of $V$. The matrix $X := A_{x\wedge y}$ satisfies
equation {\rm (\ref{form-equation-2})} of Theorem \ref{main-theorem} if and only if either $S$ is totally
isotropic or $S\cap \mathrm{Rad}(\varphi) \neq 0$.
\end{corollary}
{\bf Proof.} Since $\Phi^\sigma X^\sigma \Phi = (\Phi X \Phi^\sigma)^\sigma$ equation (\ref{form-wedge-radical})
implies (\ref{form-equation-2}). The claim immediately follows from this remark, Theorem \ref{main-theorem} and
Corollary \ref{main-corollary-bis}. \hfill $\Box$

\bigskip

Each of the matrix equations (\ref{form-equation-1}), (\ref{form-equation-2}), (\ref{form-equation-1bis}),
(\ref{form-equation-2bis}) and (\ref{form-wedge-radical}) packs in a unique formal container a number of scalar
equations. For instance, (\ref{form-equation-1bis}) is equivalent to the following set of equations:
\begin{equation}\label{form-equation-2bisbis}
\sum_{i,j=1}^nx_{k,i}^\sigma \varphi_{i,j}x_{j,h} = 0
\end{equation}
for $k, h = 1,2,...,n$, with the implicit convention that $x_{j,i} = - x_{i,j}$ and $x_{i,i} = 0$. Needless to
say, in each particular case some or even many of these $n^2$ equations become trivial or follow from other
equations of this set, possibly combined with those of (\ref{Grassmann-equazioni}), which define the Grassmann
variety $\cal G$. It can also happen that some or even all of the equations (\ref{Grassmann-equazioni}) follow
from (\ref{form-equation-2bisbis}).

\subsection{The variety ${\cal G}_\varphi$}

Let $\varphi$ be a $(\sigma,\varepsilon)$-sesquilinear form on $V$, as in the previous subsection, but now we
assume that $\varphi$ has Witt index $m\geq 2$, namely at least one line of $\mathrm{PG}(V)$ is totally
isotropic for $\varphi$, and that $\varphi$ is not the null form. Denoted by ${\cal L}_\varphi$ the set of lines
of $\mathrm{PG}(V)$ totally isotropic for $\varphi$, put ${\cal G}_\varphi := \iota_2({\cal L}_\varphi)
\subseteq {\cal G}$. By the assumptions made on $\varphi$, we have $\emptyset \neq {\cal G}_\varphi \subset
{\cal G}$. We call ${\cal G}_\varphi$ the $\varphi$-{\em subset} of $\cal G$.

When $\varphi$ is non-degenerate the set ${\cal G}_\varphi$ is described by the equations packed in
(\ref{form-equation-2bis}) of Corollary \ref{main-corollary} combined with those of (\ref{Grassmann-equazioni}).

Let $\varphi$ be degenerate. Let $\mathrm{St}(\mathrm{Rad}(\varphi)) =
\cup_{[x]\in[\mathrm{Rad}(\varphi)]}\mathrm{St}([x])$ be the set of lines of $\mathrm{PG}(V)$ meeting
$[\mathrm{Rad}(\varphi)]$ non-trivially and put ${\cal R}_\varphi :=
\iota_2(\mathrm{St}(\mathrm{Rad}(\varphi)))$. As remarked in Subsection \ref{properties grassmannian}, for every
point $[x]$ of $\mathrm{PG}(V)$ the set $\iota_2(\mathrm{St}([x]))$ is an $(n-2)$-dimensional subspace of
$\mathrm{PG}(V)$. Thus, ${\cal R}_\varphi$ is a union of subspaces of $\mathrm{PG}(V)$. Chosen a basis $U =
\{u_1,..., u_n\}$ of $V$ is such a way that $\langle u_1,..., u_r\rangle = \mathrm{Rad}(\varphi)$, the subspace
$\langle {\cal R}_\varphi\rangle$ of $V^\wedge$ is spanned by the vectors $u_i\wedge u_j$ with $i < j$ and
$i\leq r\rangle$. Hence $\mathrm{dim}(\langle{\cal R}_\varphi\rangle) = (n-r)r+{r\choose 2}$. The equations
packed in (\ref{form-wedge-radical}) decribe $\langle {\cal R}_\varphi\rangle$ while the set ${\cal R}_\varphi$
is described by those equations together with (\ref{Grassmann-equazioni}). By Corollary \ref{main-corollary-ter}
the equations packed in (\ref{form-equation-2}) combined with those of (\ref{Grassmann-equazioni}) describe
${\cal R}_\varphi\cup{\cal G}_\varphi$. We have ${\cal R}_\varphi\cup{\cal G}_\varphi = {\cal G}_\varphi$ if and
only if all points of $\mathrm{PG}(V)$ are isotropic, namely $\varphi$ is alternating. In this case
(\ref{form-equation-2}) and (\ref{Grassmann-equazioni}) characterize ${\cal G}_\varphi$.

The subgroup of $\mathrm{GL}(V)$ preserving $\varphi$ acts transitively on the set of totally
$\varphi$-isotropic lines of $\mathrm{PG}(V)$ that meet $\mathrm{Rad}(\varphi)$ trivially. Moreover, the
embedding of $\mathrm{Gr}_2(\mathrm{PG}(V))$ in $\mathrm{PG}(V^\wedge)$ is homogeneous. Therefore the setwise
stabilizer of ${\cal G}_\varphi$ in $\mathrm{PGL}(V^\wedge)$ acts transitively on ${\cal
G}_\varphi\setminus{\cal R}_\varphi$.

\subsubsection{The algebraic variety $\widetilde{\cal G}_\varphi$}\label{equation tangent section}

Let $\varphi$ be non-degenerate or alternating and suppose that either $\sigma = \mathrm{id}_{\mathbb{F}}$ or
$\mathbb{F} = \mathbb{F}_{q^2}$ for a prime power $q$ and $t^\sigma = t^q$. Then ${\cal G}_\varphi$ is the set
of $\mathbb{F}$-rational points of the (possibly reducible) algebraic variety $\widetilde{\cal G}_\varphi$
defined over the algebraic closure $\overline{\mathbb{F}}$ of $\mathbb{F}$ by the equations
(\ref{form-equation-2bisbis}) and (\ref{Grassmann-equazioni}) (if $\varphi$ is non-degenerate) or those packed
in (\ref{form-equation-2}) together with (\ref{Grassmann-equazioni}) (if $\varphi$ is alternating but
degenerate).

Put $\overline{V} = \overline{\mathbb{F}}\otimes V$ and $\overline{V}^\wedge = \overline{V}\wedge\overline{V}
\cong \overline{\mathbb{F}}\otimes V^\wedge$. Let $\overline{\cal G}\subset \mathrm{PG}(\overline{V}^\wedge)$ be
the Grassmann variety defined over $\overline{\mathbb{F}}$.

Let $\sigma = \mathrm{id}_{\mathbb{F}}$. In this case a unique bilinear form $\bar{\varphi}$ exists on
$\overline{V}$ inducing $\varphi$ on $V$, and $\widetilde{\cal G}_\varphi$ is the $\bar{\varphi}$-subset
$\overline{\cal G}_{\bar{\varphi}}$ of $\overline{\cal G}$.

Let $\mathbb{F} = \mathbb{F}_{q^2}$ and $t^\sigma = t^q$. Let $\tilde{\sigma}$ be the natural extension of
$\sigma$ to $\overline{\mathbb{F}}$, defined by the clause $t^{\tilde{\sigma}} = t^q$ for every $t\in
\overline{\mathbb{F}}$. Then $\tilde{\sigma}$ is an automorphism of $\overline{\mathbb{F}}$ but
$\tilde{\sigma}^2 \neq \mathrm{id}_{\overline{\mathbb{F}}}$. Indeed the elements of $\overline{\mathbb{F}}$
fixed by $\tilde{\sigma}^2$ are just those of $\mathbb{F}$. Consequently there is no sesquilinear form
$\bar{\varphi}$ on $\overline{V}$ such that $\widetilde{\cal G}_\varphi$ is the $\bar{\varphi}$-subset of
$\overline{\cal G}$. Nevertheless we can extend $\sigma$ to an involutory automorphism $\bar{\sigma}$ of
$\overline{\mathbb{F}}$, and we can always do it in infinitely many ways. Let $\bar{\varphi}$ be the unique
$(\bar{\sigma},\varepsilon)$-sesquilinear form on $\overline{V}$ inducing $\varphi$ on $V$. Then ${\cal
G}_\varphi$ is the set of $\mathbb{F}$-rational points of the $\bar{\varphi}$-subset $\overline{\cal
G}_{\bar{\varphi}}$ of $\overline{\cal G}$. However $\widetilde{\cal G}_\varphi \neq \overline{\cal
G}_{\bar{\varphi}}$ for any choice of $\bar{\sigma}$. Moreover, since no non-trivial involutory automorphism of
an infinite field can be expressed as a polynomial function, no polynomial $p(t)\in \overline{\mathbb{F}}[t]$
exists such that $t^{\bar{\sigma}} = p(t)$ for every $t\in\overline{\mathbb{F}}$. Thus, the set $\overline{\cal
G}_{\bar{\varphi}}$ is not an algebraic subvariety of $\overline{\cal G}$.

Finally, we warn that the variety $\widetilde{\cal G}_\varphi$ might be reducible, even disconnected, and it
might contain some singularities. We will see two examples of this kind in Section \ref{Grassmann-Emb-GQ}. In
Subsection \ref{Section 4 H(3,q)} the variety $\widetilde{\cal G}_\varphi$ is just a finite set of points while
in Subsection \ref{Section 4 Q+} it consists of two disjoint conics. In the first case $\mathbb{F} =
\mathbb{F}_{q^2}$ and $\sigma(t) = t^q$ while $\sigma = \mathrm{id}_{\mathbb{F}}$ in the second case.

The case where $\mathbb{F}$ is infinite and $\sigma \neq \mathrm{id}_{\mathbb{F}}$ remains to consider. Let
$\mathbb{F}_\sigma$ be the subfield of $\mathbb{F}$ formed by the elements fixed by $\sigma$. Since $\sigma^2 =
\mathrm{id}_{\mathbb{F}} \neq \sigma$, the field $\mathbb{F}$ is an algebraic extension of $\mathbb{F}_\sigma$
of degree $2$. Thus we can regard it as $2$-dimensional vector space over $\mathbb{F}_\sigma$ and $\sigma$ as a
linear mapping of that vector space. Accordingly, $V^\wedge$ is turned into an $n(n-1)$-dimensional vector space
$V'$ over $\mathbb{F}_\sigma$ and we can regard ${\cal G}_\varphi$ as a subset of $\mathrm{PG}(V')$. In this way
${\cal G}_\varphi$ is the set of $\mathbb{F}_\sigma$-rational points of an algebraic variety defined over the
algebraic closure $\overline{\mathbb{F}}$ of $\mathbb{F}_\sigma$ (equal to the algebraic closure of
$\mathbb{F}$, since $\mathbb{F}$ is algebraic over $\mathbb{F}_\sigma$).

\subsubsection{Tangent spaces and dimension}\label{tangent section}

Assume that $\varphi$ is either non-degenerate or alternating and either $\sigma = \mathrm{id}_{\mathbb{F}}$ or
$\mathbb{F}$ is finite. Then we can define the tangent space $\mathrm{Tan}({\cal G}_\varphi)_{[a\wedge b]}$ of
${\cal G}_\varphi$ at a point $[a\wedge b]$ of ${\cal G}_\varphi$ as the subspace of $\mathrm{PG}(V^{\wedge})$
formed by the $\mathbb{F}$-rational points of the tangent space $\mathrm{Tan}(\widetilde{\cal
G}_\varphi)_{[a\wedge b]}$ of the algebraic variety $\widetilde{\cal G}_\varphi$ at the point $[a\wedge b]$.
Thus $\mathrm{Tan}({\cal G}_\varphi)_{[a\wedge b]}$ is the subspace of $V^\wedge$ described by the linear
equations (\ref{Grassmann-equazioni-tangente}) together with those that can be obtained by differentiating
(\ref{form-equation-2bis}) (when $\varphi$ is non-degenerate) or (\ref{form-equation-2}) (if $\varphi$ is
degenerate and alternating).

Explicitly, from (\ref{form-equation-2bis}) with $\sigma = \mathrm{id}_{\mathbb{F}}$ we obtain
\begin{equation}\label{form-equation-2bis-tangent1}
X\Phi A_{a\wedge b} + A_{a\wedge b}\Phi X = O.
\end{equation}
From (\ref{form-equation-2bis}) with $\mathbb{F} = \mathbb{F}_{q^2}$ and $t^\sigma = t^q$ we get
\begin{equation}\label{form-equation-2bis-tangent2}
A_{a\wedge b}^q\Phi X = O
\end{equation}
(where if $A_{a\wedge b} = (a_{i,j})_{i,j=1}^n$ we put $A^q_{a\wedge b} = (a_{i,j}^q)_{i,j=1}^n$). Finally, from
(\ref{form-equation-2}) with $\sigma = \mathrm{id}_{\mathbb{F}}$ we get
\begin{equation}\label{form-equation-1bis-tangent}
\Phi X\Phi A_{a\wedge b} + \Phi A_{a\wedge b}\Phi X = O.
\end{equation}
By Corollary \ref{main-corollary-bis}, when $[a\wedge b]\in{\cal R}_\varphi$ the previous equation reduces to
$\Phi X\Phi A_{a\wedge b} = O$, which in its turn is implied by (\ref{form-wedge-radical}). So, if $[a\wedge
b]\in{\cal R}_\varphi$ then $\langle {\cal R}_\varphi\rangle\subseteq \mathrm{Tan}({\cal G}_\varphi)_{[a\wedge
b]}$.

If all tangent spaces of ${\cal G}_\varphi$ have the same dimension, as when $\varphi$ is non-degenerate, we
define the {\em dimension} of ${\cal G}_\varphi$ as the projective dimension of its tangent spaces.

\subsubsection{A warning}

In the equations (\ref{form-equation-1})-(\ref{form-wedge-radical}) of Subsection \ref{form-section} the unknown
$X$ stands for an arbitrary anti-symmetric matrix. Hence $X^T = -X$. Moreover $\Phi^T = \varepsilon\Phi^\sigma$
and $\sigma^2 = \mathrm{id}_{\mathbb{F}}$. It follows that $(X^\sigma\Phi X)^T = \varepsilon(X^\sigma\Phi
X)^\sigma$. Consequently, $X^\sigma\Phi X = O$ if and only if $(X^\sigma \Phi X)^T = O$. So, we might conclude
that there is no loss in taking $k \leq h$ in (\ref{form-equation-2bisbis}). This conclusion is correct if
$\sigma = \mathrm{id}_{\mathbb{F}}$ but it is wrong when $\mathbb{F} = \mathbb{F}_{q^2}$ and $t^\sigma = t^q$.
Indeed in the latter case if we differentiate the equations (\ref{form-equation-2bisbis}) to compute a system of
linear equations for $\mathrm{Tan}({\cal G}_\varphi)_{[a\wedge b]}$, but we only take those with $k \leq h$,
then we miss some of the equations that we can obtain by differentiating all of (\ref{form-equation-2bisbis}).

Here is an explanation of this puzzle. Let $\mathbb{F} = \mathbb{F}_{q^2}$ with $t^\sigma = t^q$ and let
$\tilde{\sigma}$ be the natural extension of $\sigma$ to $\overline{\mathbb{F}}$, as in Subsection \ref{equation
tangent section}. Then $\tilde{\sigma}^2\neq \mathrm{id}_{\overline{\mathbb{F}}}$. Therefore, in general
$(X^{\tilde{\sigma}}\Phi X)^T \neq \varepsilon(X^{\tilde{\sigma}}\Phi X)^{\tilde{\sigma}}$. Hence the equations
$X^{\tilde{\sigma}}\Phi X = O$ and $(X^{\tilde{\sigma}}\Phi X)^T = O$ are not equivalent over
$\overline{\mathbb{F}}$. Consequently, the subset of (\ref{form-equation-2bisbis}) formed by the equations with
$k\leq h$ describes an algebraic variety $\widehat{\cal G}_\varphi$ larger than $\widetilde{\cal G}_\varphi$.
Accordingly, if we only differentiate those equation then we obtain a description of the space
$\mathrm{Tan}(\widehat{\cal G}_\varphi)_{[a\wedge b]}$, which properly contains $\mathrm{Tan}(\widetilde{\cal
G}_\varphi)_{[a\wedge b]}$.

Of course, in principle there is nothing wrong in choosing $\widehat{\cal G}_\varphi$ instead of
$\widetilde{\cal G}_\varphi$. Indeed ${\cal G}_\varphi$ is also the set of $\mathbb{F}$-rational points of
$\widehat{\cal G}_\varphi$. But we should not forget that, when $\sigma \neq \mathrm{id}_{\mathbb{F}}$ and
$\mathbb{F}$ is finite, the tangent spaces of ${\cal G}_\varphi$ (whence the dimension of ${\cal G}_\varphi$)
depend on which of $\widetilde{\cal G}_\varphi$ or $\widehat{\cal G}_\varphi$ we choose. In Section
\ref{Grassmann-Emb-GQ}, where we will apply the theory developed so far to a number of special cases,
$\widetilde{\cal G}_\varphi$ is our choice.

\subsubsection{The span of ${\cal G}_\varphi$}

Still assuming that $\varphi$ has Witt index $m \geq 2$ and it is not the null form, in this subsection we also
assume that $\varphi$ is trace-valued, namely $V$ is spanned by the set of vectors isotropic for $\varphi$. We
recall that $\varphi$ is trace-valued if and only if $\varphi(x,x)\in\{t+\varepsilon t^\sigma\}_{t\in
\mathbb{F}}$ for every $x\in V$ (Tits \cite[Chapter 8]{Tits}). The latter condition is satisfied whenever either
$\mathrm{char}(\mathbb{F})\neq 2$ or $\sigma\neq \mathrm{id}_{\mathbb{F}}$ (recall that $\mathbb{F}$ is a field,
namely it is commutative).

\begin{theorem}\label{main-theorem-bis}
Under the previous hypotheses, let $\varphi$ be non-degenerate or alternating.

If $\varphi$ is neither degenerate nor alternating then $\langle{\cal G}_\varphi\rangle = V^\wedge$.

If $\varphi$ is alternating (possibly degenerate) then $\langle{\cal G}_\varphi\rangle$ is a hyperplane of $V^\wedge$.
\end{theorem}
{\bf Proof.} Let $\varphi$ be non-alternating, whence non-degenerate by the hypotheses of the theorem. Modulo
multiplying $\varphi$ by a suitable scalar when $\sigma \neq \mathrm{id}_{\mathbb{F}}$ we may assume that
$\varepsilon = 1$. We argue by induction on $n$, $n = 4$ being the smallest case to consider.

Suppose firstly that $n > 4$. We may assume to have chosen the basis $E$ of $V$ in such a way that
\[\varphi((x_i)_{i=1}^n, (y_i)_{i=1}^n) = \sum_{i=1}^m(x_{2i-1}^\sigma y_{2i}+ x_{2i}^\sigma y_{2i-1}) + \sum_{i=2m+1}^n\lambda_ix_i^\sigma y_i\]
for an anisotropic reflexive $\sigma$-hermitian form $\sum_{i=2m+1}^n\lambda_ix_i^\sigma y_i$ in the unknowns
$x_{2m+1},..., x_n, y_{2m+1},..., y_n$ when $n > 2m$ and with $\sum_{i=2m+1}^n\lambda_ix_i^\sigma y_i := 0$ when
$n = 2m$. Suppose firstly that $n > 2m$. Let $H_1 = \mathrm{Ker}(e^*_n)$ and $H_2 = \mathrm{Ker}(e^*_1+e^*_n)$
and let $\varphi_1$ and $\varphi_2$ be the forms induced by $\varphi$ on $H_1$ and $H_2$ respectively. Then
$\varphi_1$ and $\varphi_2$ are non-degenerate, non-alternating and trace-valued and have Witt index equal to
$m$. By the inductive hypothesis, $\langle {\cal G}_{\varphi_1}\rangle = H_1\wedge H_1$ and $\langle{\cal
G}_{\varphi_2}\rangle = H_2\wedge H_2$. We have $(H_1\wedge H_1)\cap(H_2\wedge H_2) = (H_1\cap
H_2)\wedge(H_1\cap H_2)$. It follows that
\[\mathrm{dim}((H_1\wedge H_1)+(H_2\wedge H_2)) = 2{{n-1}\choose 2}-{{n-2}\choose 2} = {n\choose 2}-1\]
namely $H^\wedge := (H_1\wedge H_1)+(H_2\wedge H_2) = \langle{\cal G}_{\varphi_1}\cup{\cal
G}_{\varphi_2}\rangle$ is a hyperplane of $V^\wedge$. However ${\cal G}_\varphi \supseteq {\cal
G}_{\varphi_1}\cup{\cal G}_{\varphi_2}$. Hence $\langle {\cal G}_\varphi\rangle\supseteq H^\wedge$. In order to
prove that $\langle{\cal G}_\varphi\rangle = V^\wedge$ it remains to show that ${\cal G}_\varphi$ also contains
a vector of $V^\wedge\setminus H^\wedge$.

The spaces $H_1\wedge H_1$ and $H_2\wedge H_2$ are represented by the equations $x_{i,n} = 0$ for $i = 1, 2,...,
n-1$ and respectively $x_{1,i}+ x_{i,n} = 0$ for $i = 1, 2,..., n-1$. Note that, since $x_{n,n} = x_{1,1} = 0$,
the equations $x_{1,1}+x_{1,n} = 0$ and $x_{1,n}+ x_{n,n} = 0$ are both equivalent to $x_{1,n} = 0$. So,
$H^\wedge$ is described by the equation $x_{1,n} = 0$. Put now $u = e_1$ and $v = e_{2m-1} + e_{2m}s + e_n$ for
a scalar $s$ such that $s+s^\sigma + \lambda_n = 0$. A scalar $s$ with this property exists because $\lambda_n =
\varphi(e_n,e_n)\in\{t+t^\sigma\}_{t\in\mathbb{F}}$ (recall that $\varphi$ is trace-valued and $\varepsilon =
1$). If $t\in \mathbb{F}$ is such that $\lambda_n = t + \varepsilon t^\sigma$ then take $s = -t$.

With this choice of $u$ and $v$ the $2$-space $\langle u, v\rangle$ is totally isotropic. Hence $[u\wedge v]\in
{\cal G}_\varphi$. On the other hand, $u\wedge v = e_{1,2m-1}+ e_{1,2m}s + e_{1,n}$, which does not satisfy the
equation $x_{1,n} = 0$ of $H^\wedge$. Hence $u\wedge v\not \in H^\wedge$, as we wanted.

The case where $n = 2m$ can be dealt with in a similar way, considering the hyperplanes $H_1 =
\mathrm{Ker}(e^*_1+e^*_2)$ and $H_2 = \mathrm{Ker}(e^*_{2m-1}+e^*_m)$. The forms $\varphi_1$ and $\varphi_2$
induced on $H_1$ and $H_2$ have Witt index $m-1$, but now $m > 2$ since $2m = n > 4$. Hence $m-1 \geq 2$ and the
induction hypothesis can be applied. The hyperplane $H^\wedge = (H_1\wedge H_1)+(H_2\wedge H_2)$ is described by
the equation $x_{1,n-1}+ x_{1,n} + x_{2,n-1} +x_{2,n} = 0$. If we take $u = e_1$ and $v = e_n$ then $\langle u,
v\rangle$ is totally isotropic but $u\wedge v$ does not satisfies the equation of $H^\wedge$.

Finally, let $n = 2m = 4$. In this case we must find six totally isotropic 2-subspaces $\langle u_i, v_i\rangle$
for $i = 1, 2,..., 6$ such that the vectors $u_i\wedge v_i$ are linearly independent. The following choice does
the job:
\[\begin{array}{ll}
\begin{array}{lr}
u_1 = e_1, & v_1 = e_3,\\
u_3 = e_1, & v_3 = e_4,
\end{array} &
\begin{array}{lr}
u_2 = e_2, & v_2 = e_4,\\
u_4 = e_2, & v_4 = e_3,
\end{array} \\
~ u_5 = e_1-e_3, & ~ v_5 = e_2+ e_4, \\
u_6 = -e_1t+e_2+ e_3t+e_4, & v_6 = -e_1t+ e_3(t-t^\sigma) + e_4
\end{array}\]
where $t+t^\sigma = 1$. Note that in the present case if $\mathrm{char}(\mathbb{F}) = 2$ then $\sigma \neq
\mathrm{id}_{\mathbb{F}}$ otherwise, since $\varepsilon = 1$ and $n = 2m = 4$, the form $\varphi$ is
alternating, contrary to our assumptions. So, the equation $t+t^\sigma = 1$ admits a solution.

Assume now that $\varphi$ is alternating. If $\varphi$ is non-degenerate then it is well known that ${\cal
G}_\varphi$ spans a hyperplane of $V^\wedge$ (see Blok and Cooperstein \cite{BCoop}, for instance). We have
assumed that $m \geq 2$ but the previous claim holds true for $m = 1$ too. Indeed in this case ${\cal G}_\varphi
= \emptyset$ but $V^\wedge$ is $1$-dimensional. Suppose now that $\varphi$ is degenerate. Then
\[\varphi((x_i)_{i=1}^n,(y_i)_{i=1}^n) = \sum_{i=1}^k(x_{2i-1}y_{2i}-x_{2i}y_{2i-1})\]
for a positive integer $k$ such that $2k < n$ and $\mathrm{Rad}(\varphi) = \langle e_{2k+1},..., e_n\rangle$.
Let $S = \langle e_1,..., e_{2k}\rangle$. Then $\varphi$ induces a non-degenerate alternating form $\varphi_S$
on $S$. By the above, ${\cal G}_{\varphi_S}$ spans a hyperplane of $S^\wedge := S\wedge S$. On the other hand,
${\cal G}_\varphi$ contains ${\cal R}_\varphi$, which is a subspace of $\mathrm{PG}(V^\wedge)$. Regarded ${\cal
R}_\varphi$ as a subspace of $V^\wedge$, as we may, we have $V^\wedge = S^\wedge\oplus{\cal R}_\varphi$. So,
$\langle {\cal G}_\varphi\rangle = \langle {\cal G}_{\varphi_S}\rangle \oplus{\cal R}_\varphi$ is a hyperplane
of $V^\wedge$.  \hfill $\Box$

\subsubsection{The set $\iota_2(\mathrm{St}([x]))\cap{\cal G}_\varphi$
for a point $[x]\in{\cal S}_\varphi$}\label{immagine St(x)}

In this subsection we recall a few well known facts, for further reference. Denoted by ${\cal S}_\varphi$ the
set of points of $\mathrm{PG}(V)$ isotropic for $\varphi$ let $\overline{\cal S}_{\varphi} = {\cal
S}_\varphi\setminus[\mathrm{Rad}(\varphi)]$ ($= {\cal S}_\varphi$ when $\varphi$ is non-degenerate). The group
of linear transformations of $V$ that preserve $\varphi$ acts transitively on $\overline{\cal S}_\varphi$ and on
$[\mathrm{Rad}(\varphi)]$ (when $\varphi$ is degenerate). So, chosen a point $[a]\in\overline{\cal S}_\varphi$
and a point $[b]\in[\mathrm{Rad}(\varphi)]$, the set $\iota_2(\mathrm{St}([x]))\cap{\cal G}$ is isomorphic to
either $\iota_2(\mathrm{St}([a]))\cap{\cal G}$ or $\iota_2(\mathrm{St}([b]))\cap{\cal G}$, according to whether
$[x]\in \overline{\cal S}_\varphi$ or $[x]\in[\mathrm{Rad}(\varphi)]$.

We may assume to have chosen the basis $E = (e_1,..., e_n)$ of $V$ in such a way that $[e_n]\in{\cal
S}_\varphi$. If $e_n\in \mathrm{Rad}(\varphi)$ then we can take $b = e_n$, otherwise $a = e_n$. With this choice
of $E$ the matrix $\Phi$ can be described as follows
\[\Phi = \left[\begin{array}{cc}
\Phi_0 & \upsilon \\
\varepsilon(\upsilon^\sigma)^T & 0
\end{array}\right]\]
where $\Phi_0$ is an $(n-1)\times(n-1)$ matrix such that $\Phi_0^T = \varepsilon \Phi_0^\sigma$ and $\upsilon$
is an $(n-1)\times 1$ matrix. We have $\upsilon = 0$ if and only if $e_n\in \mathrm{Rad}(\varphi)$. Given a
non-zero vector $x = \sum_{i=1}^ne_ix_i$, we have
\[A_{x\wedge e_n} = \left[\begin{array}{cc}
O & \hat{x} \\
-\hat{x}^T & 0
\end{array}\right]\]
where $O$ stands for the $(n-1)\times(n-1)$ null matrix and $\hat{x} = (x_i)_{i=1}^{n-1}$, written as an
$(n-1)\times 1$ matrix. If $e_n\not\in\mathrm{Rad}(\varphi)$ (whence $\upsilon\neq 0$) by equation
(\ref{form-equation-2}) with $X = A_{x\wedge e_n}$ we obtain that $\langle x,e_n\rangle$ is totally isotropic if
and only if
\begin{equation}\label{form-equation-star}
(\upsilon^\sigma)^T\hat{x} = 0 ~~~~ \mbox{and} ~~~~ (\hat{x}^\sigma)^T\Phi_0\hat{x} = 0.
\end{equation}
The set $\iota_2(\mathrm{St}([e_n]))\cap{\cal G}_\varphi$ is described by the equations
(\ref{form-equation-star}) together with the linear equations (\ref{Grassmann-equazioni-star}) describing the
subspace $S_{[e_n]} = \iota_2(\mathrm{St}([e_n]))$ of $\mathrm{PG}(V^\wedge)$. If $\varphi$ is trace-valued then
$\iota_2(\mathrm{St}([e_n])\cap{\cal G}_\varphi$ spans the hyperplane $S_{[e_n]}^\upsilon$ of $S_{[e_n]}$
described by the first equation of (\ref{form-equation-star}). When $\varphi$ is alternating the second equation
of (\ref{form-equation-star}) is trivial. In this case $\iota_2(\mathrm{St}([e_n]))\cap{\cal G}_\varphi =
S_{[e_n]}^\upsilon$.

Let $e_n\in\mathrm{Rad}(\varphi)$. Then (\ref{form-equation-2}) trivially holds for any $X = A_{x\wedge e_n}$,
hence it gives no conditions on $\hat{x}$. However in this case $\langle x, e_n\rangle$ is totally isotropic if
and only if $x-e_nx_n$ is isotropic. As the entries of $\hat{x}$ are just the first $n-1$ coordinates of
$x-e_nx_n$, the vector $x-e_nx_n$ is isotropic if and only if the second equation of (\ref{form-equation-star})
holds. The first one is trivial, since now $\upsilon = 0$.

If $\varphi$ is alternating and $e_n\in\mathrm{Rad}(\varphi)$ then both equations (\ref{form-equation-star}) are
trivial.

\subsection{Lines totally singular for a quadratic form}

Let $\chi$ be a pseudoquadratic form on $V$ and $\varphi$ its sesquilinearization (Tits \cite[Chapter 8]{Tits}).
Since $\mathbb{F}$ is a field, if either $\mathrm{char}(\mathbb{F}) \neq 2$ or $\sigma \neq
\mathrm{id}_{\mathbb{F}}$ then we have $\chi(x) = 0$ if and only if $\varphi(x,x) = 0$. In this case there is no
need to consider $\chi$ in addition to $\varphi$. Thus, in this subsection we assume that
$\mathrm{char}(\mathbb{F}) = 2$ and $\sigma = \mathrm{id}_{\mathbb{F}}$, whence $\varepsilon = 1$. So $\chi$ is
a quadratic form. We also assume that $\chi$ is non-singular and has Witt index $m \geq 2$, namely at least one
line of $\mathrm{PG}(V)$ is totally singular for $\chi$. The bilinearization $\varphi$ of $\chi$ is a possibly
degenerate alternating form and $\chi(x)\neq 0$ for every non-zero vector $x\in\mathrm{Rad}(\varphi)$, if any.
Let $\Phi$ be the representative matrix of $\varphi$.

\begin{theorem}\label{main-theorem-quadric}
For two non-proportional vectors $x = \sum_{i=1}^nx_i$ and $y = \sum_{i=1}^ne_iy_i$ of $V$, let $X =
(x_{i,j})_{i,j=1}^n := A_{x\wedge y}$. The $2$-subspace $\langle x,y\rangle$ of $V$ is totally singular for
$\chi$ if and only $\Phi X\Phi X = O$ and
\begin{equation}\label{form-equation-quadric}
\chi((x_{i,k})_{i=1}^n) = 0 ~~~ \mbox{for}~ $k = 1, 2,..., n$.
\end{equation}
\end{theorem}
{\bf Proof.} Recalling that $x_{i,j} = x_iy_j + x_jy_i$, it is easy to see that
\begin{equation}\label{quadric equation}
\chi((x_{i,k})_{i=1}^n) = y_k^2\chi(x) + x_k^2\chi(y) + x_ky_k\varphi(x,y).
\end{equation}
Put $S = \langle x,y\rangle$. By Corollary \ref{main-corollary-ter} and since ${\cal R}_\varphi\subseteq{\cal
G}_\varphi$ because $\varphi$ is alternating, the space $S$ is totally isotropic for $\varphi$ if and only if
$\Phi X \Phi X = O$. It is totally singular for $\chi$ if and only if it is totally isotropic and $\chi(x) =
\chi(y) = 0$, if and only if $\chi(x) = \chi(y) = \varphi(x,y) = 0$. So, if $S$ is totally singular it follows
from (\ref{quadric equation}) that $\chi((x_{i,k})_{i=1}^n) = 0$ for every $k = 1, 2,..., n$. Conversely,
assuming that $S$ is totally isotropic for $\varphi$, suppose that all equations (\ref{form-equation-quadric})
hold for $S$. Modulo replacing $x$ and $y$ with other two non-proportional vectors of $S$ if necessary, we may
assume that two indices $r < s$ exist such that $x_r = 1$, $x_s = 0$, $y_r = 0$ and $y_s = 1$. Then $x_{i,r} =
y_i$ and $x_{i,s} = x_i$ for $i = 1, 2,..., n$. By (\ref{form-equation-quadric}) and (\ref{quadric equation})
with $k = r$ and $k = s$ we obtain $\chi(y) = 0$ and $\chi(x) = 0$, respectively. So, $S$ is totally singular.
\hfill $\Box$

\bigskip

Note that in general the equation $\Phi X\Phi X = O$ is not a consequence of (\ref{form-equation-quadric}) and
(\ref{Grassmann-equazioni}), as we can see by the following example.

\begin{example}\label{Ex}
{\rm Let $n = 4$ and $\chi(x_1, x_2, x_3, x_4) = x_1x_2 + x_3x_4$. Then the equations $\Phi X\Phi X = O$,
(\ref{form-equation-quadric}) and (\ref{Grassmann-equazioni}) can be written as follows:
\[\begin{array}{ll}
\left.\begin{array}{ll}
(A.1) & x_{1,2}^2 + x_{1,4}x_{2,3} + x_{1,3}x_{2,4} = 0,\\
(A.2) & x_{1,2} = x_{3,4},
\end{array}\right\} & (\mbox{from}~ \Phi X\Phi X = O)\\
\left.\begin{array}{ll}
(B.1) & x_{1,3}x_{1,4} = 0,\\
(B.2) & x_{2,3}x_{2,4} = 0,\\
(B.3) & x_{1,3}x_{2,3} = 0,\\
(B.4) & x_{1,4}x_{2,4} = 0,
\end{array}\right\} & (\mbox{from}~ (\ref{form-equation-quadric}))\\
\begin{array}{ll}
~~ (C) & x_{1,2}x_{3,4} + x_{1,3}x_{2,4} + x_{1,4}x_{2,3} = 0.
\end{array} & (\mbox{from}~ (\ref{Grassmann-equazioni}))
\end{array}\]
Modulo $(A.2)$, the equations $(C)$ and $(A.1)$ are equivalent, but there is no way to deduce $(A.2)$ from
$(B.1), (B.2), (B.3), (B.4)$ and $(C)$.}
\end{example}

\subsubsection{The variety ${\cal G}_\chi$}

Denoted by ${\cal L}_\chi$ the set of lines of $\mathrm{PG}(V)$ totally singular for $\chi$, put ${\cal G}_\chi
:= \iota_2({\cal L}_\chi) \subseteq {\cal G}_\varphi$. By Theorem \ref{main-theorem-quadric} the set ${\cal
G}_\chi$ is a projective variety described by the equations (\ref{form-equation-quadric}) together with those
packed in $\Phi X \Phi X = O$ and the equations (\ref{Grassmann-equazioni}) describing $\cal G$. Its tangent
space at a point $[a\wedge b]$ is described by the linear equations (\ref{Grassmann-equazioni-tangente})
together with those packed in $\Phi A_{a\wedge b} \Phi X + \Phi X \Phi A_{a\wedge b} = O$ (compare
(\ref{form-equation-1bis-tangent})) and the following ones, obtained by differentiating
(\ref{form-equation-quadric}):
\begin{equation}\label{form-equation-quadric-tan}
\sum_{j=1}^n\left[\frac{\partial \chi((x_{i,k})_{i=1}^n)}{\partial x_{j,k}}\right]_{a\wedge b}x_{j,k} = 0 ~~~ \mbox{for}~ k = 1, 2,..., n.
\end{equation}
Clearly, $\langle {\cal G}_\chi\rangle$ is contained in the hyperplane $\langle {\cal G}_\varphi\rangle$ of
$V^\wedge$ (compare Theorem \ref{main-theorem-bis}).

\begin{theorem}\label{span of Gq}
We have $\langle{\cal G}_\chi\rangle = \langle{\cal G}_\varphi\rangle$.
\end{theorem}
{\bf Proof.} The proof is quite similar to that of Theorem \ref{main-theorem-bis}. We will give only a sketch of
it, leaving details for the reader. We may assume to have chosen the basis $E$ of $V$ in such a way that
\[\chi((x_i)_{i=1}^n) = \sum_{i=1}^m x_{2i-1}x_{2i} + \chi_0(x_{2m+1},...., x_n)\]
where if $n > 2m$ then $\chi_0$ is a totally anisotropic quadratic form as follows
\[\chi_0(x_{2m+1},..., x_n) = \sum_{i=m+1}^r(\lambda_{2i-1}x_{2i-1}^2 + x_{2i-1}x_{2i}+ \lambda_{2i}x_{2i}^2) + \sum_{i= 2r+1}^n \lambda_ix_i^2\]
while $\chi_0 := 0$ if $n = 2m$. Note that if $\mathbb{F}$ is perfect then $n-2m \leq 2$, but in the general
case $n$ can be larger than $m+2$. We allow $r = m$ and $n = 2r$.

We argue by induction on $n$. Suppose firstly that $n > 4$. Let $n = 2m$. Put $H_1 = \mathrm{Ker}(e^*_1+e^*_2)$
and $H_2 = \mathrm{Ker}(e^*_{n-1}+e^*_n)$ and $H^\wedge = (H_1\wedge H_1)+(H_2\wedge H_2)$. By the inductive
hypothesis on the forms induced by $\chi$ on $H_1$ and $H_2$, the intersection $H^\wedge\cap\langle{\cal
G}_\varphi\rangle$ is a hyperplane of $\langle{\cal G}_\varphi\rangle$. In order to prove that ${\cal G}_\chi$
spans $\langle{\cal G}_\varphi\rangle$ we only must find a point of ${\cal G}_\chi$ outside $H^\wedge$. The
point $[e_1\wedge e_n]$ has the required properties.

Suppose now that $n > 2r \geq 2m$. In this case we choose $H_1 = \mathrm{Ker}(e^*_n)$ and $H_2 =
\mathrm{Ker}(e_1^*+e_n^*)$. With $H^\wedge$ defined as above, let $v = e_{2m-1}+e_{2m}\lambda_n + e_n$. Then
$[e_1\wedge v]$ is a point of ${\cal G}_\chi$ outside $H^\wedge$.

Finally, let $n = 2r > 2m$. Put $H_1 = \mathrm{Ker}(e^*_{n-1})$, $H_2 = \mathrm{Ker}(e^*_n)$ and $v = e_{2m-1} +
e_{2m}(\lambda_{n-1}+\lambda_n + 1) + e_{n-1}+ e_n$. The point $[e_1\wedge v]$ belongs to ${\cal G}_\chi$ but
not to $H^\wedge$.

The case of $n = 4$ remains to examine. In this case ${\cal G}_\chi$ contains the points represented by the
vectors $e_{1,3}$, $e_{1,4}$, $e_{2,3}$, $e_{2,4}$ and $e_{1,2}+e_{1,3}+e_{2,4}+e_{3,4} =
(e_2+e_3)\wedge(e_1+e_4)$ (compare the equations describing ${\cal G}_\chi$ in Example \ref{Ex}). These five
vectors are linearly independent. So, they span $\langle{\cal G}_\varphi\rangle$, which is $5$-dimensional.
\hfill $\Box$

\bigskip

Theorem \ref{span of Gq} answers in the affirmative Conjectures 4.5 and 4.12 of \cite{Pas-survey} (see also
\cite[Conjecture 1]{CP1}) in the special case of $k = 2$.

\subsubsection{The set $\iota_2(\mathrm{St}([a]))\cap{\cal G}_\chi$ for a point $[a]\in{\cal G}_\chi$}\label{residuo in car 2}

The group of linear transformations of $V$ preserving $\chi$ acts transitively on ${\cal G}_\chi$. Therefore,
chosen a point $[a]\in{\cal G}_\chi$, we have $\iota_2(\mathrm{St}([x]))\cap{\cal
G}_\chi\cong\iota_2(\mathrm{St}([a]))\cap{\cal G}_\chi$ for every point $[x]\in{\cal G}_\chi$. Assuming to have
chosen the basis $E$ of $V$ in such a way that $[e_n]\in{\cal G}_\chi$, take $a = e_n$. As in Subsection
\ref{immagine St(x)},
\[\Phi = \left[\begin{array}{cc}
\Phi_0 & \upsilon \\
\varepsilon(\upsilon^\sigma)^T & 0
\end{array}\right]\]
for an $(n-1)\times 1$ matrix $\upsilon$. Note that $\upsilon\neq 0$ since $[e_n]\in{\cal G}_\chi$ and ${\cal
G}_\chi\cap[\mathrm{Rad}(\varphi)] = \emptyset$. Given a vector $x = \sum_{i=1}e_ix_i$ of $V$ non-proportional
to $e_n$ let $\hat{x} := (x_i)_{i=1}^{n-1}$ and put $\chi_0(\hat{x}) := \chi(x_1,...,x_{n-1},0)$. The 2-subspace
$\langle x,e_n\rangle$ is totally singular if and only if
\begin{equation}\label{form-equation-star-quadric}
(\upsilon^\sigma)^T\hat{x} = 0 ~~~~ \mbox{and} ~~~~ \chi_0(\hat{x}) = 0.
\end{equation}

\subsection{Further investigations}

Most likely, the results of this section can be generalized to $k$-grassmannians of polar spaces of rank $m \geq
k$. The vectors of $\wedge^kV$ correspond to alternating tensors of degree $k$. The crucial step in this project
is to find a way to generalize Lemma \ref{fundamental lemma} to $k$-subspaces of $V$ and $V^*$.

\section{Grassmann embeddings of quadrangles}\label{Grassmann-Emb-GQ}

Throughout this section $\Gamma = ({\cal P}, {\cal L})$ is a generalized quadrangle, $\Gamma^* = ({\cal P}^*,
{\cal L}^*)$ is its dual (see Subsection \ref{Grassmann embeddings}), $\varepsilon$ is a projective embedding of
$\Gamma^*$ in $\mathrm{PG}(V)$ for a vector space $V$ of dimension $n\geq 4$ over a given field $\mathbb{F}$ and
$\varepsilon^{\mathrm{gr}}$ is the Grassmann embedding of $\Gamma$ induced by $\varepsilon$ (see Subsection
\ref{Grassmann embeddings}). We keep the notation of Section \ref{Section 3}. In particular, $V^\wedge :=
V\wedge V$ and, given an ordered basis $E = (e_1,..., e_n)$ of $V$, the vectors $e_{i,j} = e_i\wedge e_j$ for $i
< j$ form the basis $E^\wedge$ of $V^\wedge$ canonically associated to $E$. Moreover, $\iota_2$ is the natural
embedding of $\mathrm{Gr}_2(\mathrm{PG}(V))$ in $\mathrm{PG}(V^\wedge)$ and ${\cal G} =
\iota_2(\mathrm{Gr}_2(\mathrm{PG}(V)))$.

Given a trace-valued non-degenerate $(\sigma, \varepsilon)$-sesquilinear form $\varphi$ on $V$ of Witt index
$2$, we denote by $Q(\varphi)$ the generalized quadrangle associated to $\varphi$, namely $Q(\varphi) = ({\cal
P}_\varphi, {\cal L}_\varphi)$ where ${\cal P}_\varphi$ is the set of points of $\mathrm{PG}(V)$ isotropic for
$\varphi$ and ${\cal L}_\varphi$ is the set of lines of $\mathrm{PG}(V)$ totally isotropic for $\varphi$.
Similarly, given a non-singular quadratic form $\chi$ of Witt index $2$, $Q(\chi) = ({\cal P}_\chi,{\cal
L}_\chi)$ is the generalized quadrangle formed by the points and the lines of $\mathrm{PG}(V)$ that are singular
and respectively totally singular for $\chi$. Also, ${\cal G}_\varphi := \iota_2({\cal L}_\varphi)$ and ${\cal
G}_\chi := \iota_2({\cal L}_\chi)$, as in Section \ref{Section 3}. So, if $\varepsilon(\Gamma^*) = Q(\varphi)$
then $\varepsilon^{\mathrm{gr}}$ embeds $\Gamma$ in $\langle {\cal G}_\varphi\rangle^{\mathrm{pr}}$ and we have
$\varepsilon^{\mathrm{gr}}({\cal P}) = {\cal G}_\varphi$. We recall that $\langle{\cal
G}_\varphi\rangle^{\mathrm{pr}}$ is a hyperplane of $\mathrm{PG}(V^\wedge)$ when $\varphi$ is alternating while
$\langle{\cal G}_\varphi\rangle^{\mathrm{pr}} = \mathrm{PG}(V^\wedge)$ in all other cases (Theorem
\ref{main-theorem-bis}).

Similarly, if $\mathrm{char}(\mathbb{F}) = 2$ and $\varepsilon(\Gamma^*) = Q(\chi)$ for a quadratic form $\chi$
then $\varepsilon^{\mathrm{gr}}$ embeds $\Gamma$ in $\langle {\cal G}_\chi\rangle^{\mathrm{pr}}$ and we have
$\varepsilon^{\mathrm{gr}}({\cal P}) = {\cal G}_\chi$, but now $\langle{\cal G}_\chi\rangle^{\mathrm{pr}}$ is a
hyperplane of $\mathrm{PG}(V^\wedge)$ (Theorem \ref{span of Gq}).

We shall examine several special cases of the situation described above. In each of the cases we will consider
but the last two $\Gamma$ also admits a projective embedding $\varepsilon_0$. We will also consider the
veronesean embedding $\varepsilon_0^{\mathrm{ver}}$ induced by $\varepsilon_0$ (Subsection \ref{veronesean
embeddings}), discussing its relations with $\varepsilon^{\mathrm{gr}}$, if any.

\subsection{Quadrangles of symplectic type}\label{Section 4 W(3,q)}

Let $\Gamma \cong W(3,\mathbb{F})$, the generalized quadrangle associated to a non-degenerate alternating form
of $V(4,\mathbb{F})$. Then $\Gamma^* \cong Q(4,\mathbb{F})$. Accordingly, $V$ has dimension $n = 5$ and
$\varepsilon(\Gamma^*) = Q(\chi)$ with $\chi$ a non-singular quadratic form of $V$ of Witt index $2$. Let
$\varphi$ be the bilinearization of $\chi$. If $\mathrm{char}(\mathbb{F})\neq 2$ then $Q(\chi) = Q(\varphi)$. By
(\ref{form-equation-star}) of Subsection \ref{immagine St(x)} when $\mathrm{char}(\mathbb{F})\neq 2$ and
(\ref{form-equation-star-quadric}) of Subsection \ref{residuo in car 2} when $\mathrm{char}(\mathbb{F}) = 2$ it
follows that $\varepsilon^{\mathrm{gr}}(l)$ is a non-singular conic for every line $l$ of $\Gamma$. Thus, the
embedding $\varepsilon^{\mathrm{gr}}$ is quadratic.

\subsubsection{The case $\mathrm{char}(\mathbb{F}) \neq 2$}

Let $\mathrm{char}(\mathbb{F})\neq 2$. Then $Q(\chi) = Q(\varphi)$. The form $\varphi$ is symmetric. By Theorem
\ref{main-theorem-bis}, the set $\varepsilon^{\mathrm{gr}}({\cal P}) = {\cal G}_\varphi$ ($= {\cal G}_\chi$)
spans $V^\wedge \cong V(10,\mathbb{F})$.

Modulo multiplying $\varphi$ by a scalar, we can assume to have chosen the basis $E$ of $V$ in such a way that
$\varphi$ is represented by the following matrix with respect to $E$:
\[\Phi = \left[\begin{array}{ccccc}
0 & 1 & 0 & 0 & 0 \\
1 & 0 & 0 & 0 & 0 \\
0 & 0 & 0 & 1 & 0 \\
0 & 0 & 1 & 0 & 0 \\
0 & 0 & 0 & 0 & 1
\end{array}\right].\]
With the convention that $x_{i,i} = 0$ and $x_{j,i} = - x_{i,j}$, the variety ${\cal G}_\varphi$ is described by
the equations (\ref{Grassmann-equazioni}) together with the following ones, packed in $X\Phi X = O$ (Corollary
\ref{main-corollary}, (\ref{form-equation-2bis})):
\[x_{k,1}x_{2,h} + x_{k,2}x_{1,h} + x_{k,3}x_{4,h}+x_{k,4}x_{3,h} + x_{k,5}x_{5,h} = 0, ~~~ (1\leq k \leq h \leq 5).\]
The equation $A_{a\wedge b}\Phi X + X\Phi A_{a\wedge b} = O$, together with
(\ref{Grassmann-equazioni-tangente}), describes the tangent space $\mathrm{Tan}({\cal G}_\varphi)_{[a\wedge b]}$
of ${\cal G}_\varphi$ at a point $[a\wedge b]\in{\cal G}_\varphi$. Clearly, $\varepsilon^{\mathrm{gr}}$ is
homogeneous. So, we are free to choose the point $[a\wedge b]$ as we like. We choose $a = e_1$ and $b = e_3$, so
that $a\wedge b = e_{1,3}$. With this choice of $[a\wedge b]$ we get the following linear system of rank 6 for
$\mathrm{Tan}({\cal G}_\varphi)_{[a\wedge b]}$:
\[x_{1,2} - x_{3,4} = x_{1,4} = x_{2,3} = x_{2,4} = x_{2,5} = x_{4,5} = 0.\]
As $\mathrm{dim}(V^\wedge) = 10$ we obtain that $\mathrm{dim}(\mathrm{Tan}({\cal G}_\varphi)_{[a\wedge b]}) = 4$
(vector dimension). Consequently, ${\cal G}_\varphi$ is $3$-dimensional.

\subsubsection{The case $\mathrm{char}(\mathbb{F}) = 2$}

Let $\mathrm{char}(\mathbb{F}) = 2$. By Theorems \ref{main-theorem-bis} and \ref{span of Gq}, the set
$\varepsilon^{\mathrm{gr}}({\cal P}) = {\cal G}_\chi$ spans a hyperplane of $V^\wedge$.

Modulo multiplying $\chi$ by a scalar we can assume to have chosen the basis $E$ of $V$ in such a way that
$\chi(x_1,..., x_5) = x_1x_2 + x_3x_4 + x_5^2$. Accordingly, the bilinearization $\varphi$ of $\chi$ is
represented by the following matrix:\[\Phi = \left[\begin{array}{ccccc}
 0 & 1 & 0 & 0 & 0 \\
1 & 0 & 0 & 0 & 0 \\
0 & 0 & 0 & 1 & 0 \\
0 & 0 & 1 & 0 & 0 \\
0 & 0 & 0 & 0 & 0
\end{array}\right].\]
The variety ${\cal G}_\chi$ is described by the following
equations, obtained from $\Phi X\Phi X = O$ (Corollary
\ref{main-corollary-ter}) and (\ref{form-equation-quadric}):
\[\begin{array}{ll}
x_{k,1}x_{2,h} + x_{k,2}x_{1,h} + x_{k,3}x_{4,h}+x_{k,4}x_{3,h} = 0, & (1\leq k < h\leq 5);\\
x_{1,k}x_{2,k} + x_{3,k}x_{4,k} + x_{5,k}^2 = 0, & (k = 1,...,5).
\end{array}\]
It is not difficult to check that these equations imply (\ref{Grassmann-equazioni}). Moreover, they are
equivalent to the following set of equations:
\[(A) ~~~~  x_{1,2} + x_{3,4} = 0,\]
\[(B) ~ \left\{\begin{array}{l}
x_{1,2}^2 + x_{1,3}x_{2,4} + x_{1,4}x_{2,3} = 0,\\
x_{1,5}^2 + x_{1,3}x_{1,4} = 0, \\
x_{2,5}^2 + x_{2,3}x_{2,4} = 0, \\
x_{3,5}^2 + x_{1,3}x_{2,3} = 0, \\
x_{4,5}^2 + x_{1,4}x_{2,4} = 0,
\end{array}\right.\]
\[(C) ~ \left\{\begin{array}{l}
x_{1,5}x_{2,5} + x_{3,5}x_{4,5} = 0,\\
x_{1,2}x_{1,5} + x_{1,3}x_{4,5}+ x_{1,4}x_{3,5} = 0,\\
x_{1,2}x_{2,5} + x_{2,3}x_{4,5}+ x_{2,4}x_{3,5} = 0,\\
x_{1,3}x_{2,5} + x_{2,3}x_{1,5}+ x_{1,2}x_{3,5} = 0,\\
x_{1,4}x_{2,5} + x_{2,4}x_{1,5}+ x_{1,2}x_{4,5} = 0.
\end{array}\right.\]
Equation $(A)$ describes the hyperplane of $V^\wedge$ spanned by ${\cal G}_\chi$. Note that the squares of the
equations listed in $(C)$ follow from the equations $(B)$, but this fact does not allow us to drop $(C)$
(compare Hilbert's Nullstellensatz). In order to compute $\mathrm{Tan}({\cal G}_\chi)_{[a\wedge b]}$ we choose
$a = e_1$ and $b = e_3$, as in the case of $\mathrm{char}(\mathbb{F})\neq 2$. We get just the same linear system
as in that case, namely:
\[x_{1,2} - x_{3,4} = x_{1,4} = x_{2,3} = x_{2,4} = x_{2,5} = x_{4,5} = 0.\]
(Note that if we drop the set of equations $(C)$ then we miss the last two equations of this linear system.)
So, ${\cal G}_\chi$ has dimension $3$ in this case too.

\subsubsection{The veronesean embedding $\varepsilon_0^{\mathrm{ver}}$}

As $\Gamma \cong W(3,\mathbb{F})$, a projective embedding
$\varepsilon_0:\Gamma\rightarrow\mathrm{PG}(3,\mathbb{F})$ also exists such that $\varepsilon_0(\Gamma) =
Q(\alpha)$ for a non-degenerate alternating form $\alpha$ of $V(4,\mathbb{F})$. Since $\varepsilon_0({\cal P})$
is equal to the set of points of $\mathrm{PG}(3,\mathbb{F})$, the image $\varepsilon_0^{\mathrm{ver}}({\cal P})$
of $\cal P$ by the veronesean embedding $\varepsilon_0^{\mathrm{ver}}$ induced by $\varepsilon_0$ is equal to
the Veronese variety ${\cal V} = \eta^{\mathrm{ver}}(\mathrm{PG}(3,\mathbb{F}))$.

As proved in \cite{CP1} (see also \cite{CP2} for more details), we have $\varepsilon^{\mathrm{ver}}_0\geq
\varepsilon^{\mathrm{gr}}$, namely there exists a linear mapping $\pi$ from $\langle {\cal V}\rangle =
V(10,\mathbb{F})$ to $\langle\varepsilon^{\mathrm{gr}}({\cal P})\rangle$ inducing a bijection from $\cal V$ to
$\varepsilon^{\mathrm{gr}}({\cal P})$ and mapping every conic $\varepsilon_0^{\mathrm{ver}}(l)$ for $l\in{\cal
L}$ onto the corresponding conic $\varepsilon^{\mathrm{gr}}(l)$. We call $\pi$ the {\em projection} of
$\varepsilon_0^{\mathrm{ver}}$ onto $\varepsilon^{\mathrm{gr}}$.

It is proved in \cite{CP1} that $\pi$ is an isomorphism if and only if $\mathrm{char}(\mathbb{F}) \neq 2$ and
that $\mathrm{Ker}(\pi)$ is $1$-dimensional when $\mathbb{F}$ is a perfect field of characteristic $2$. In
\cite{CP2} it is shown that the nuclei of the conics $\varepsilon_0^{\mathrm{ver}}(l)$ for $l\in{\cal L}$ form a
copy $\cal N$ of the quadric $Q(4,\mathbb{F})$ ($\cong \Gamma^*$) and $\mathrm{Ker}(\pi)$ contains the nucleus
$n({\cal N})$ of $\cal N$. So, if $\mathbb{F}$ is perfect then $\mathrm{Ker}(\pi) = n({\cal N})$.

Theorems \ref{main-theorem-bis} and \ref{span of Gq} of the present paper allows us to improve the above result.
Indeed, by those two theorems we immediately obtain the following:

\begin{proposition}\label{projection pi}
If $\mathrm{char}(\mathbb{F})\neq 2$ then $\pi$ is an isomorphism. If $\mathrm{char}(\mathbb{F}) = 2$ then
$\mathrm{Ker}(\pi)$ is $1$-dimensional, whence $\mathrm{Ker}(\pi) = n({\cal N})$.
\end{proposition}
Consequently, and since $\varepsilon_0^{\mathrm{ver}}({\cal P}) = {\cal V}$,

\begin{corollary}\label{propjection pi bis}
If $\mathrm{char}(\mathbb{F})\neq 2$ then the Veronese variety $\cal V$ of $\mathrm{PG}(3,\mathbb{F})$ is
isomorphic to the subvariety ${\cal G}_\varphi$ of $\cal G$.
\end{corollary}
Suppose now that $\mathrm{char}(\mathbb{F}) = 2$. By Proposition \ref{projection pi}, the arguments used in
\cite{CP1} and \cite{CP2} to study $\pi$ when $\mathbb{F}$ is perfect work in general, thus yielding the
following:

\begin{corollary}\label{projection pi ter}
Let $\mathrm{char}(\mathbb{F}) = 2$. Then the image $\pi({\cal N})$ of $\cal N$ by $\pi$ is isomorphic to
$W(3,\mathbb{F})$ and consists of the nuclei of the conics $\varepsilon^{\mathrm{gr}}(l)$ for $l\in{\cal L}$.
The vector space $K := \langle \pi({\cal N})\rangle$ ($ = \pi(\langle {\cal N}\rangle)$) is $4$-dimensional and
defines a quotient of $\varepsilon^{\mathrm{gr}}$. The quotient $\varepsilon^{\mathrm{gr}}/K$ is a $1$-embedding
of $\Gamma$ in $\mathrm{PG}(4,\mathbb{F})$.
\end{corollary}
The image of $\Gamma$ by $\varepsilon^{\mathrm{gr}}/K$ is a subgeometry of a copy of $\Gamma^*$ in
$\mathrm{PG}(4,\mathbb{F})$ (compare De Bruyn and Pasini \cite[Theorem 3.11]{DB-Pas}. We recall that if
$\mathbb{F}$ is perfect then $\Gamma\cong \Gamma^*$, whence $\varepsilon$ is also a projective embedding of
$\Gamma$. In this case $\varepsilon^{\mathrm{gr}}/K \cong \varepsilon$ (whence $\varepsilon^{\mathrm{gr}}/K$
maps $\Gamma$ onto a copy of $\Gamma^*$) and the restriction of $\pi$ to $\langle{\cal N}\rangle$ is the natural
projection of $\varepsilon$ onto $\varepsilon_0$.

In any case, $\pi$ induces a bijective morphism $\pi_{|{\cal V}}$ from the Veronese variety $\cal V$ to ${\cal
G}_\chi$. We do not know if $\pi_{|{\cal V}}$ is an isomorphism of algebraic varieties when
$\mathrm{char}(\mathbb{F}) = 2$. We guess it is not.

Turning to universality, we mention the following result, proved in \cite{CP2}:

\begin{proposition}\label{W(3,q) univ}
Let $\mathbb{F} = \mathbb{F}_q$ for a prime power $q > 3$. Then $\varepsilon_0^{\mathrm{ver}}$ ($\cong
\varepsilon^{\mathrm{gr}}$ by Proposition \ref{projection pi}) is relatively universal.
\end{proposition}
We guess that this statement holds true for any field of characteristic different from $2$. On the other hand,
if $\mathrm{char}(\mathbb{F}) = 2$ then $\varepsilon_0^{\mathrm{ver}}$ is not relatively universal. Indeed in
this case $\Gamma$ also admits a projective embedding $\varepsilon_{\delta}$ in
$\mathrm{PG}(4+\delta,\mathbb{F})$ for any $\delta\leq [\mathbb{F}:\mathbb{F}^2]$, where $\mathbb{F}^2$ stands
for the subfield of square elements of $\mathbb{F}$ (De Bruyn and Pasini \cite{DB-Pas}). Clearly, if
$\mathbb{F}$ is perfect then $\delta = 1$, otherwise $\delta$ can be larger than 1, even infinite. To make
things easier, when $[\mathbb{F}:\mathbb{F}^2]$ is infinite assume to have chosen a finite value for $\delta$.
The image of $\cal P$ by $\varepsilon_\delta$ is a non-singular quadric $Q_\delta$ and $\varepsilon_0\cong
\varepsilon_\delta/N_\delta$, where $N_\delta$ is the radical of the bilinearization of $Q_\delta$, namely
$[N_\delta]$ is the nucleus of $Q_\delta$. We can always consider the veronesean embedding
$\varepsilon_\delta^{\mathrm{ver}}$ induced by $\varepsilon_\delta$. The projection from $\varepsilon_\delta$ to
$\varepsilon_0$ lifts to a morphism from $\varepsilon_\delta^{\mathrm{ver}}$ to $\varepsilon_0^{\mathrm{ver}}$,
but $\mathrm{dim}(\varepsilon_\delta^{\mathrm{ver}}) > \mathrm{dim}(\varepsilon_0^{\mathrm{ver}})$. Hence
$\varepsilon_0^{\mathrm{ver}}$ cannot be relatively universal.

In particular, if $[\mathbb{F}:\mathbb{F}^2]$ is infinite then the hull of $\varepsilon_0^{\mathrm{ver}}$ is
infinite dimensional. This is not surprising, since in this case the hull of $\varepsilon_0$ is infinite
dimensional as well.

\subsection{Quadrics in $\mathrm{PG}(4,\mathbb{F})$}

In this subsection we consider a setting dual to that of Subsection \ref{Section 4 W(3,q)}. Now $\Gamma \cong
Q(4,\mathbb{F})$, the generalized quadrangle associated to a non-singular quadratic form of $V(5,\mathbb{F})$ of
Witt index $2$. Accordingly, $\Gamma^* \cong W(3,\mathbb{F})$ and $\varepsilon$ is a projective embeddings of
$W(3,\mathbb{F})$.

We recall that when $\mathrm{char}(\mathbb{F})\neq 2$ the natural embedding of $W(3,\mathbb{F})$ in
$\mathrm{PG}(3,\mathbb{F})$ (uniquely determined up to isomorphisms) is the unique projective embedding of
$W(3,\mathbb{F})$. As noticed at the end of the previous subsection, this is no more true when
$\mathrm{char}(\mathbb{F}) = 2$. In this case we must say which is the embedding $\varepsilon$ of $\Gamma^*$
that we consider. We assume that $\varepsilon$ is the natural one in this case too. So, in either case $V$ is
$4$-dimensional and $\varepsilon(\Gamma^*) = Q(\varphi)$ for a non-degenerate alternating form $\varphi$.

\subsubsection{The embedding $\varepsilon^{\mathrm{gr}}$}

It is well known that $\varepsilon^{\mathrm{gr}}$ is isomorphic to the natural projective embedding
$\varepsilon_0$ of $\Gamma$ as $Q(4,\mathbb{F})$ in $\mathrm{PG}(4,\mathbb{F})$. This conclusion can be obtained
from the results of Section \ref{Section 3} of this paper, as we are going to show.

The space $V^\wedge$ is $6$-dimensional and ${\cal G}_\varphi$ spans a hyperplane of $V^\wedge$, by Theorem
\ref{main-theorem-bis}. So, $\mathrm{dim}(\langle {\cal G}_\varphi\rangle) = 5$ (vector dimension). The equation
$\hat{x}^T\Phi_0\hat{x} = 0$ of (\ref{form-equation-star}) of Subsection \ref{immagine St(x)} is trivial, since
$\Phi_0$ is anti-symmetric. Consequently, $\varepsilon^{\mathrm{gr}}(l)$ is a projective line. Thus,
$\varepsilon^{\mathrm{gr}}$ is projective. Chosen the representative matrix $\Phi$ of $\varphi$ as follows,
\[\Phi = \left[\begin{array}{rrrr}
 0 & 0 & 1 & 0  \\
0 & 0 & 0 & 1 \\
-1 & 0 & 0 & 0 \\
0 & -1 & 0 & 0
\end{array}\right]\]
the matrix equation $X\Phi X = O$ of Corollary \ref{main-corollary} is equivalent to the following pair of
scalar equations:
\[\begin{array}{l}
x_{1,3}+ x_{2,4} = 0 ~~~ (\mbox{which describes}~ \langle{\cal G}_\varphi\rangle),\\
x_{1,4}x_{2,3} + x_{1,2}x_{2,4} + x_{1,3}^2 = 0.
\end{array}\]
Hence $\varepsilon^{\mathrm{gr}}\cong \varepsilon_0$. Note that all equations (\ref{Grassmann-equazioni}) follow
from the previous two equations.

\subsubsection{The veronesean embedding $\varepsilon_0^{\mathrm{ver}}$}

The quadric veronesean embedding $\eta^{\mathrm{ver}}$ embeds $\mathrm{PG}(4,\mathbb{F})$ in
$\mathrm{PG}(14,\mathbb{F})$ and $\varepsilon_0^{\mathrm{ver}}({\cal P})$ spans a hyperplane $H$ of
$\mathrm{PG}(14,\mathbb{F})$. Explicitly, let $\varepsilon_0({\cal P})$ be the quadric described by the equation
$x_1x_2 + x_3x_4 + x_5^2 = 0$, as we can always assume. Then $H$ is described by the equation $x_{1,2}+x_{3,4}+
x_{5,5} = 0$.

When $\mathrm{char}(\mathbb{F}) \neq 2$ no morphism can exist between $\varepsilon_0^{\mathrm{ver}}$ and
$\varepsilon_0$ ($\cong \varepsilon^{\mathrm{gr}}$). Let $\mathrm{char}(\mathbb{F}) = 2$ and let $N$ be the
subspace of $\mathrm{PG}(14,\mathbb{F})$ spanned by the nuclei of the conics $\eta^{\mathrm{ver}}(l)$ for $l$ a
line of $\mathrm{PG}(4,\mathbb{F})$. Let $N_\Gamma$ be the subspace of $N$ spanned by the nuclei of the conics
$\varepsilon_0^{\mathrm{gr}}(l) = \eta^{\mathrm{ver}}(\varepsilon_0(l))$ for $l$ a line of $\Gamma$. Clearly,
$N_\Gamma\subseteq H\cap N$.

It is known that $N$ has codimension $5$ in $\mathrm{PG}(14,\mathbb{F})$, it defines a quotient of
$\eta^{\mathrm{ver}}$ and the quotient $\eta^{\mathrm{ver}}/N$ is a $1$-embedding of $\mathrm{PG}(4,\mathbb{F})$
in itself (Thas and Van Maldeghem \cite{TVM2004}). In particular, if $\mathbb{F}$ is perfect then
$\eta^{\mathrm{ver}}/N$ is isomorphic to the identity embedding of $\mathrm{PG}(4,\mathbb{F})$ in itself.
Moreover, the mapping $\nu$ that maps every line of $\mathrm{PG}(4,\mathbb{F})$ onto its nucleus is isomorphic
to the projective embedding $\iota_2$ of $\mathrm{Gr}_2(\mathrm{PG}(4,\mathbb{F}))$ in
$\mathrm{PG}(9,\mathbb{F})$ ($\cong N$).

According to the above, the subspace $N_H := H\cap N$ defines a quotient of $\varepsilon^{\mathrm{ver}}_0$.
Clearly, $\varepsilon^{\mathrm{ver}}_0/N_H$ is a $1$-embedding of $\Gamma$.

\begin{lemma}\label{lemma NH}
The subspace $N_H$ is a hyperplane of $N$.
\end{lemma}
{\bf Proof.} We must prove that $N\not\subseteq H$, namely there exists a line $l$ of
$\mathrm{PG}(4,\mathbb{F})$ such that $\nu(l)\not\in H$. Let $l$ be the line described by the linear system $x_3
= x_4 = x_5 = 0$. Then $\nu(l)$ is the point $[x_{i,j}]_{1\leq i\leq j\leq 5}$ with $x_{i,j} = 0$ for every
$(i,j) \neq (1,2)$. This point does not satisfies the equation of $H$.  \hfill $\Box$.

\bigskip

By Lemma \ref{lemma NH}, the subspace $N_H\subset H$ has codimension $5$ in $H$. The claims gathered in the next
proposition follow from this fact, the remarks before Lemma \ref{lemma NH} and those after Corollary
\ref{projection pi ter} (modulo permuting the roles of $\varepsilon$ and $\varepsilon_0$). We leave the details
of the proof to the reader.

\begin{proposition}
We have $N_\Gamma = N_H$, $\mathrm{dim}(\varepsilon^{\mathrm{ver}}_0/N_H) = 4$ and the composition
$\nu\cdot\varepsilon_0^{\mathrm{ver}}$ is isomorphic to the Grassmann embedding $\varepsilon_0^{\mathrm{gr}}$ of
$\Gamma^* \cong W(3,\mathbb{F})$ induced by $\varepsilon_0$ (recall that $\Gamma$ is the dual of $\Gamma^*$).
Moreover, if $\mathbb{F}$ is perfect then $\varepsilon_0^{\mathrm{ver}}/N_H\cong \varepsilon_0$.
\end{proposition}
Let $\mathbb{F}$ be non-perfect. Then the $1$-embedding $\varepsilon_0^{\mathrm{ver}}/N_H$ is lax and the image
of $\Gamma$ by $\varepsilon_0^{\mathrm{ver}}/N_H$ is a proper subgeometry of $\varepsilon_0(\Gamma)$. Hence
$\Gamma$ is isomorphic to a proper subgeometry of itself. Clearly, the same holds for $\Gamma^*$. Moreover, each
of the geometries $\Gamma$ and $\Gamma^*$ is isomorphic to a proper subgeometry of the other one (De Bruyn and
Pasini \cite{DB-Pas}).

\subsection{Hermitian varieties in $\mathrm{PG}(3,\mathbb{K})$}

Given a field $\mathbb{F}$, let $\mathbb{K}$ be an extension of $\mathbb{F}$ of degree $[\mathbb{K}:\mathbb{F}]
= 2$ and $\sigma$ the unique involutory automorphism of $\mathbb{K}$ fixing $\mathbb{F}$ elementwise. Let
$\Gamma$ be isomorphic to the generalized quadrangle associated to a non-degenerate $\sigma$-hermitian form of
$V(4,\mathbb{K})$ of Witt index $2$. It is well known that $\Gamma^* \cong Q^-(5,\mathbb{F})$. So, $V$ is a
$6$-dimensional vector space over $\mathbb{F}$ and $\varepsilon(\Gamma^*) = Q(\chi)$ for a non-singular
quadratic form $\chi$ of Witt index $2$ with non-degenerate bilinearization $\varphi$. Recall that if
$\mathrm{char}(\mathbb{F}) \neq 2$ then ${\cal G}_\chi = {\cal G}_\varphi$ while ${\cal G}_\chi \subset {\cal
G}_\varphi$ when $\mathrm{char}(\mathbb{F}) \neq 2$.

The space $V^\wedge$ is $15$-dimensional. By Theorems \ref{main-theorem-bis} and \ref{span of Gq}, if
$\mathrm{char}(\mathbb{F}) \neq 2$ then ${\cal G}_\varphi$ spans $V^\wedge$ while $\langle {\cal G}_\chi\rangle
= \langle{\cal G}_\varphi\rangle$ is a hyperplane of $V^\wedge$ when $\mathrm{char}(\mathbb{F}) = 2$.

By (\ref{form-equation-star}) of Subsection \ref{immagine St(x)} when $\mathrm{char}(\mathbb{F}) \neq 2$ and
(\ref{form-equation-star-quadric}) of Subsection \ref{residuo in car 2} when $\mathrm{char}(\mathbb{F}) = 2$ it
follows that $\mathrm{dim}(\langle\varepsilon^{\mathrm{gr}}(l)\rangle^{\mathrm{pr}}) = 3$ for every line $l$ of
$\Gamma$ and $\varepsilon^{\mathrm{gr}}(l)$ is an elliptic quadric of the $3$-space
$\langle\varepsilon^{\mathrm{gr}}(l)\rangle^{\mathrm{pr}}$. Thus $\varepsilon^{\mathrm{gr}}$ is a $3$-embedding.

Let $\mathrm{char}(\mathbb{F}) \neq 2$. We can assume that $\varphi$ is represented by the following matrix, for
a suitable non-square element $\eta$ of $\mathbb{F}$:
\[\Phi = \left[\begin{array}{rrrrrr}
 0 & 1 & 0 & 0 & 0 & 0 \\
1 & 0 & 0 & 0 & 0 & 0 \\
0 & 0 & 0 & 1 & 0 & 0 \\
0 & 0 & 1 & 0 & 0 & 0 \\
0 & 0 & 0 & 0 & 1 & 0 \\
0 & 0 & 0 & 0 & 0 &-\eta
\end{array}\right]\]
The variety ${\cal G}_\varphi$ is described by the equations (\ref{Grassmann-equazioni}) together with the
following ones, packed in $X\Phi X = O$ (i.e. (Corollary \ref{main-corollary}, (\ref{form-equation-2bis})):
\[x_{k,1}x_{2,h} + x_{k,2}x_{1,h} + x_{k,3}x_{4,h}+x_{k,4}x_{3,h} + x_{k,5}x_{5,h} -\eta x_{k,6}x_{6,h} = 0\]
for $1\leq k\leq h\leq 6$. The equation $A_{a\wedge b}\Phi X + X\Phi A_{a\wedge b} = O$, together with
(\ref{Grassmann-equazioni-tangente}), describes the tangent space $\mathrm{Tan}({\cal G}_\varphi)_{[a\wedge
b]}$. We choose $a = e_1$ and $b = e_3$. With this choice of $[a\wedge b]$ we get the following linear system of
rank $9$ for $\mathrm{Tan}({\cal G}_\varphi)_{[a\wedge b]}$:
\[x_{1,2} - x_{3,4} = x_{1,4} = x_{2,3} = x_{2,4} = x_{2,5} = x_{2,6} = x_{4,5} = x_{4,6} = x_{5,6} = 0.\]
As $\mathrm{dim}(V^\wedge) = 15$ we obtain that $\mathrm{dim}(\mathrm{Tan}({\cal G}_\varphi)_{[a\wedge b]}) = 6$
(vector dimension). Consequently, ${\cal G}_\varphi$ is $5$-dimensional.

Let $\mathrm{char}(\mathbb{F}) = 2$. We can assume that
\[\chi(x_1,..., x_6) = x_1x_2 + x_3x_4 + x_5x_6 + x_5^2 + \lambda x_6^2\]
for a scalar $\lambda \in \mathbb{F}$ such that the polynomial $t^2 + t +\lambda$ is irreducible over
$\mathbb{F}$. Accordingly, $\varphi$ is represented by the following matrix
\[\Phi = \left[\begin{array}{rrrrrr}
 0 & 1 & 0 & 0 & 0 & 0 \\
1 & 0 & 0 & 0 & 0 & 0 \\
0 & 0 & 0 & 1 & 0 & 0 \\
0 & 0 & 1 & 0 & 0 & 0 \\
0 & 0 & 0 & 0 & 0 & 1 \\
0 & 0 & 0 & 0 & 1 & 0
\end{array}\right]\]
The variety ${\cal G}_\chi$ is described by the equations (\ref{Grassmann-equazioni}) together with the
following equations, obtained from $X\Phi X = O$ (Corollary \ref{main-corollary}) and
(\ref{form-equation-quadric}):
\[x_{k,1}x_{2,h} + x_{k,2}x_{1,h} + x_{k,3}x_{4,h}+x_{k,4}x_{3,h} + x_{k,5}x_{6,h} + x_{k,6}x_{5,h} = 0\]
for $1\leq k < h\leq 6$ and
\[x_{1,k}x_{2,k} + x_{3,k}x_{4,k} + x_{5,k}x_{6,k} + x_{5,k}^2 + \lambda x_{5,k}^2 = 0\]
for $k = 1,...,6$. The tangent space $\mathrm{Tan}({\cal G}_\chi)_{[e_{1,3}]}$ is described by the same system
of linear equations as when $\mathrm{char}(\mathbb{F}) \neq 2$. Hence ${\cal G}_\chi$ has dimension $5$.

Apparently, no relation exists between $\varepsilon_0^{\mathrm{ver}}$ and $\varepsilon^{\mathrm{gr}}$.

\subsection{Elliptic quadrics in $\mathrm{PG}(5,q)$}\label{Section 4 H(3,q)}

Let $\Gamma \cong Q^-(5,q)$, the generalized quadrangle associated to a non-singular quadratic form of $V(6,q)$
of Witt index $2$, for a prime power $q$. Then $\Gamma^* \cong H(3,q^2)$, namely $V$ is a $4$-dimensional vector
space over $\mathbb{F}_{q^2}$ and $\varepsilon(\Gamma^*) = Q(\varphi)$ for a non-degenerate hermitian form
$\varphi$ of $V$ of Witt index $2$. It is well known that $\varepsilon^{\mathrm{gr}}$ is a lax $1$-embedding. It
embeds $\Gamma$ as $Q^-(5,q)$ in a Baer subgeometry of $\mathrm{PG}(5,q^2)$ (see Cooperstein and Shult
\cite{CS-DPS}, for instance). Thus, $\varepsilon^{\mathrm{gr}}$ is obtained from the natural embedding
$\varepsilon_0:\Gamma\rightarrow\mathrm{PG}(5,q)$ by extending the field $\mathbb{F}_q$ to $\mathbb{F}_{q^2}$.

Some of the previous statements can be obtained from the results of Section \ref{Section 3} of our paper. By
Theorem \ref{main-theorem-bis} the set ${\cal G}_\varphi$ spans $V^\wedge$, which is $6$-dimensional. Given a
point $[a]\in{\cal G}_\varphi$, the space $\mathrm{St}([a])$ is a projective plane. By the first equation of
(\ref{form-equation-star}) of Subsection \ref{immagine St(x)}, $\langle
\varepsilon^{\mathrm{gr}}(l)\rangle^{\mathrm{pr}}$ is a line. By the second equation of
(\ref{form-equation-star}), $\varepsilon^{\mathrm{gr}}(l)$ is a Baer subline of $\langle
\varepsilon^{\mathrm{gr}}(l)\rangle^{\mathrm{pr}}$. Thus, $\varepsilon^{\mathrm{gr}}$ is a lax $1$-embedding.

We may assume to have chosen the basis $E$ of $V$ in such a way that $\varphi$ is represented by the identity
matrix. In this way ${\cal G}_\varphi$ is represented by the matrix equation $X^qX = O$ (Corollary
\ref{main-corollary}), combined with equations (\ref{Grassmann-equazioni}). Explicitly, $X^qX = O$ can be
written as follows:
\[x_{k,1}^qx_{1,h} + x_{k,2}^qx_{2,h} + x_{k,3}^qx_{3,h} + x_{k,4}^qx_{4,h} = 0, ~~~ (1\leq k, h \leq 4).\]
The tangent space $\mathrm{Tan}({\cal G}_\varphi)_{[a\wedge b]}$ is described by
(\ref{Grassmann-equazioni-tangente}) together with the equation $A_{a\wedge b}^qX = O$ (see
(\ref{form-equation-2bis-tangent2})). Choose $a = e_1+ e_2t$ and $b = e_3+e_4t$ for a given $t$ such that
$t^{q+1}+1 = 0$, and $t \neq 1$ if $q$ is even. So, $a\wedge b = e_{1,2} + (e_{1,4}+q_{2,3})t + e_{2,4}t^2$.
With this choice of $a\wedge b$ we get the following linear system of rank $5$ for the tangent space
$\mathrm{Tan}({\cal G}_\varphi)_{[a\wedge b]}$:
\[x_{2,3} + t^qx_{2,4} = x_{1,2} = x_{1,3} = x_{1,4} = x_{3,4} = 0.\]
As $\mathrm{dim}(V^\wedge) = 6$ we obtain that $\mathrm{dim}(\mathrm{Tan}({\cal G}_\varphi)_{[a\wedge b]}) = 1$
(vector dimension). Consequently, ${\cal G}_\varphi$ is $0$-dimensional.

In fact, let $\overline{\mathbb{F}}$ be the algebraic closure of $\mathbb{F}_{q^2}$ and put $\overline{V}^\wedge
= \overline{\mathbb{F}}\otimes V^\wedge$. It is straightforward to check that the solutions of the equation
$X^qX = O$ in $\mathrm{PG}(\overline{V}^\wedge)$ are just the same as in $\mathrm{PG}(V^\wedge)$. In other
words, $\widetilde{\cal G}_\varphi = {\cal G}_\varphi$ (notation as in Subsection \ref{equation tangent
section}). Thus, $\widetilde{\cal G}_\varphi$ is a finite set of points.

As for $\varepsilon_0^{\mathrm{ver}}$, apparently no relation exists between $\varepsilon^{\mathrm{ver}}_0$ and
$\varepsilon^{\mathrm{gr}}$.

\subsection{The dual of $H(4,q^2)$}

Let $\Gamma^*\cong H(4,q^2)$, namely $V$ is a $5$-dimensional vector space over $\mathbb{F}_{q^2}$ and
$\varepsilon(\Gamma^*) = Q(\varphi)$ for a non-degenerate hermitian form $\varphi$ of Witt index $2$. The space
$V^\wedge$ is $10$-dimensional and ${\cal G}_\varphi$ spans $V^\wedge$, by Theorem \ref{main-theorem-bis}. Given
a point $[a]\in{\cal G}_\varphi$, $\mathrm{St}([a])$ is a $3$-dimensional subspace of $\mathrm{PG}(V^\wedge)$.
By (\ref{form-equation-star}) of Subsection \ref{immagine St(x)}, the set $\varepsilon^{\mathrm{gr}}(l)$ is a
unital, spanning a plane of $\mathrm{St}([a])$. So, $\varepsilon^{\mathrm{gr}}$ is a $2$-embedding, but it is
not quadratic.

As in the previous subsection, we may assume to have chosen the basis $E$ of $V$ in such a way that $\varphi$ is
represented by the identity matrix. In this way ${\cal G}_\varphi$ is represented by the matrix equation $X^qX =
O$ combined with (\ref{Grassmann-equazioni}). Explicitly, $X^qX = O$ can be written as follows:
\[x_{k,1}^qx_{1,h} + x_{k,2}^qx_{2,h} + x_{k,3}^qx_{2,h} + x_{k,4}^qx_{4,h} + x_{k,5}^qx_{5,h} = 0, ~~~ (1\leq k, h \leq 5).\]
The tangent space $\mathrm{Tan}({\cal G}_\varphi)_{[a\wedge b]}$ is described by the matrix equation $A_{a\wedge
b}^qX = O$ together with (\ref{Grassmann-equazioni-tangente}). Choose $a = e_1+ e_2t$ and $b = e_3+e_4t$ for a
given $t$ such that $t^{q+1}+1 = 0$, and $t \neq 1$ if $q$ is even. So, $a\wedge b = e_{1,2} +
(e_{1,4}+q_{2,3})t + e_{2,4}t^2$. With this choice of $a\wedge b$ we get the following linear system of rank $7$
for $\mathrm{Tan}({\cal G}_\varphi)_{[a\wedge b]}$:
\[\begin{array}{l}
x_{1,2} = x_{1,3} = x_{1,4} = x_{3,4} = 0,\\
x_{2,3} + t^qx_{2,4} = x_{1,5} + t^qx_{2,5} = x_{3,5} + t^qx_{4,5} = 0.
\end{array}\]
As $\mathrm{dim}(V^\wedge) = 10$, it follows that $\mathrm{dim}(\mathrm{Tan}({\cal G}_\varphi)_{[a\wedge b]}) =
3$. Consequently, ${\cal G}_\varphi$ is $2$-dimensional.

It is well known that the dual of $H(4,q^2)$ admits no projective embedding. So, no projective embedding exists
for $\Gamma$.

\subsection{Dual grids}\label{Section 4 Q+}

Let $\Gamma^* \cong Q^+(3,\mathbb{F})$ and let $\varepsilon:\Gamma^*\rightarrow\mathrm{PG}(3,\mathbb{F})$ be a
projective embedding of $\Gamma^*$. Thus, $\varepsilon(\Gamma^*) = Q(\chi)$ for a non-singular quadratic form
$\chi$ of Witt index $2$ in $V = V(4,\mathbb{F})$. The space $V^\wedge$ is $6$-dimensional.

We may assume that $\chi(x_1, x_2, x_3, x_4) = x_1x_2 + x_3x_4$. We have already given the equations of ${\cal
G}_\chi$ when $\mathrm{char}(\mathbb{F}) = 2$ in Example \ref{Ex}. When $\mathrm{char}(\mathbb{F}) \neq 2$, we
still obtain the equations $(B.1)-(B.4)$ of Example \ref{Ex} but $(A.1)$ and $(A.2)$ are replaced by $x_{1,2}^2
= x_{1,3}x_{2,4} + x_{1,4}x_{2,3}$ and $x_{1,2}^2 = x_{2,4}^2$ respectively. In addition to these equations we
also obtain the following ones:
\[\begin{array}{lcl}
x_{1,3}(x_{1,2}-x_{3,4}) = 0, & & x_{2,4}(x_{1,2}-x_{3,4}) = 0, \\
x_{1,4}(x_{1,2}+x_{3,4}) = 0, & & x_{2,3}(x_{1,2}+x_{3,4}) = 0.
\end{array}\]
Moreover, $x_{1,2}x_{3,4} - x_{1,3}x_{2,4} + x_{1,4}x_{3,4} = 0$ by (\ref{Grassmann-equazioni}).

It follows that ${\cal G}_\chi$ is the union of two disjoint conics, say $C_1$ and $C_2$. When
$\mathrm{char}(\mathbb{F}) \neq 2$ the planes $\langle C_1\rangle^{\mathrm{pr}}$ and $\langle
C_2\rangle^{\mathrm{pr}}$ are disjoint. In this case $C_1\cup C_2$ spans $\mathrm{PG}(V^\wedge)$. If
$\mathrm{char}(\mathbb{F}) = 2$ then $C_1$ and $C_2$ have the same nucleus, say $n$, and $\langle
C_1\rangle^{\mathrm{pr}}\cap\langle C_2\rangle^{\mathrm{pr}} = n$. In this case $C_1\cup C_2$ spans a hyperplane
of $\mathrm{PG}(V^\wedge)$ (compare Theorems \ref{main-theorem-bis} and \ref{span of Gq}).

\bigskip

\noindent
Ilaria Cardinali and Antonio Pasini\\
Department of Information Engineering and Mathematics\\
University of Siena\\
Via Roma 56, 53100 Siena\\
cardinali3@unisi.it, pasini@unisi.it


\begin{thebibliography}{XYZ}
\bibitem{TVMmore} Z. Akca, A. Bayar, S. Ekmekci, R. Kaya, J. A. Thas and H. Van Maldeghem. Generalized veronesean embeddings of projective spaces, Part II. The lax case. {\em Ars Combin.} {\bf 103} (2012), 65-80.
\bibitem{BCoop} R. J. Blok and B. N. Cooperstein. The generating rank of the unitary and symplectic grassmannians. {\em J. Combin. Th. A} {\bf 119} (2012), 1-13.
\bibitem{BP} R. J. Blok and A. Pasini. On absolutely universal embeddings. {\em Discrete Math.} {\bf 267} (2003), 45--62.
\bibitem{CP1} I. Cardinali and A. Pasini. Grassmann and Weyl embeddings of orthogonal grassmannians. To appear in {\em J. Alg. Combin.}
\bibitem{CP2} I. Cardinali and A. Pasini. Veronesean embeddings of dual polar spaces of orthogonal type. To appear in {\em J. Combin. Th. A}.
\bibitem{CS-DPS} B. N. Cooperstein and E. E. Shult. A note on embedding and generating dual polar spaces. {Adv. Geometry} {\bf 1} (2001), 37-48.
%\bibitem{BDBpseudo} B. De Bruyn. Pseudo-embeddings and pseudo-hyperplanes. To appear in {\em Advances in Geometry}.
\bibitem{DB-Pas} B. De Bruyn and A. Pasini. On symplectic polar spaces over non-perfect fields of characteristic $2$. {\em Linear and Multilinear Algebra} {\bf 57} (2009), 567-575.
\bibitem{HT} J. W. P. Hirschfeld and J. A. Thas. {\em General Galois Geometries.} Oxford Univ. Press, Oxford 1991.
\bibitem{KS} A. Kasikova and E. E. Shult. Absolute embeddings of point-line geometries. {\em J. Algebra} {\bf 238} (2001), 265--291.
%\bibitem{PasDG} A. Pasini. {\em Diagram Geometries}. Oxford University Press, Oxford, 1994.
\bibitem{PasEE} A. Pasini. Embeddings and expansions. {\em Bull. Belg. Math. Soc. Simon Stevin} {\bf 10} (2003),  585--626.
\bibitem{Pas-survey} A. Pasini. Embeddings of orthogonal grassmannians. To appear in {\em Innovations in Incidence Geometry}.
\bibitem{PasVM} A. Pasini and H. Van Maldeghem. Some constructions and embeddings of the tilde geometry. {\em Note di Matematica} {\bf 21} (2002/2003), 1-33.
%\bibitem{PT} S. E. Payne and J. A. Thas. {\em Finite Generalized Quadrangles}. Pitman, Boston, 1984.
\bibitem{Ronan} M. A. Ronan. Embeddings and hyperplanes of discrete geometries. {\em European J. Combin.} {\bf 8}  (1987), 179--185.
\bibitem{TVM2004} J. A. Thas and H. Van Maldeghem. Characterizations of the finite quadric veroneseans ${{\mathcal{V}}_n}^{2^{n}}$.  {\em Q. J. Math.} {\bf 55} (2004), no. 1, 99–-113.
\bibitem{TVM} J. A. Thas and H. Van Maldeghem. Generalized veronesean embeddings of projective spaces. {\em Combinatorica} {\bf 31} (2011), 615--629.
\bibitem{Tits} J. Tits. {\em Buildings of Spherical Type and Finite $BN$-pairs.} Lect. Notes in Math. {\bf 386}, Springer, Berlin, 1974.
\bibitem{HVM} H. Van Maldeghem. {\em Generalized Polygons}. Birkh\"{a}user, Basel, 1998.
\end{thebibliography}
\end{document}